\documentclass[11pt,a4paper]{article}
\usepackage[round,authoryear,comma]{natbib}
\usepackage{optimstyle}
\usepackage{tabularx,pifont,calc,fontawesome5,placeins}
\newcommand{\R}{\mathbb R}
\newcommand{\E}{\mathbb E}
\newcommand{\ip}[2]{\langle #1,#2\rangle}
\newcommand{\norm}[1]{\lVert#1\rVert}
\newcommand{\method}[1]{\textcolor{optimGreen}{\texttt{#1}}}
\newcommand{\one}{\mathbf 1}
\newcommand{\Gap}{\mathtt{Gap}}
\DeclareMathOperator{\dist}{dist}

\DeclareMathOperator{\bl}{bl}
\DeclareMathOperator{\pos}{pos}
\newcommand{\by}[2]{\mathrel{\underset{#2}{#1}}}

\newcolumntype{P}[1]{>{\raggedright\arraybackslash}p{#1}}
\newcolumntype{Y}{>{\raggedright\arraybackslash}X}

\setlist[itemize]{label=$\vartriangleright$}
\algrenewcommand\algorithmicrequire{\textbf{Input:}}
\algrenewcommand\algorithmicensure{\textbf{Output:}}
\algrenewcommand\algorithmicfor{\textbf{for}}
\algrenewcommand\algorithmicwhile{\textbf{while}}
\algrenewcommand\algorithmicloop{\textbf{loop}}
\algrenewcommand\algorithmicdo{\textbf{do}}
\algrenewcommand\algorithmicend{\textbf{end}}
\algrenewcommand\algorithmicif{\textbf{if}}
\algrenewcommand\algorithmicthen{\textbf{then}}
\algrenewcommand\algorithmicelse{\textbf{else}}
\algrenewcommand\algorithmicreturn{\textbf{return}}
\BeforeBeginEnvironment{theorem}{\Needspace{6\baselineskip}}
\BeforeBeginEnvironment{lemma}{\Needspace{6\baselineskip}}
\BeforeBeginEnvironment{proposition}{\Needspace{6\baselineskip}}
\BeforeBeginEnvironment{definition}{\Needspace{6\baselineskip}}
\BeforeBeginEnvironment{assumption}{\Needspace{6\baselineskip}}
\pretocmd{\subsection}{\Needspace{8\baselineskip}}{}{}
\renewcommand{\optimtitlefont}{\fontsize{16}{21}\selectfont\bfseries\boldmath}
\title{Descending Everest, Part III\texorpdfstring{\\[3pt]}{: }
Oracle Complexity of Monotone\texorpdfstring{\\}{ }
Variational Inequalities\texorpdfstring{\\}{ }
with \texorpdfstring{$(L_0,L_1)$}{(L0,L1)}-Lipschitz Operators\texorpdfstring{\thanks{The Russian version of this paper is available at \url{https://github.com/chitoraphelion/descending-everest-with-first-order-oracles}.}}{}}
\author{%
  Mikhail Rudomanov\thanks{MIRIAI, Moscow Independent Research Institute of Artificial Intelligence, Moscow, Russia. \email{rudomanov.ma@phystech.edu}.}\,\thanks{Equal contribution.}\quad
  Ivan Fomin\thanks{Yandex, Moscow, Russia. \email{ash-ivanfomin@yandex-team.ru}.}\,\footnotemark[3]\quad
  Alexander V. Gasnikov\thanks{Innopolis University, Innopolis, Russia. \email{gasnikov@yandex.ru}.}\quad
  Fedor S. Stonyakin\thanks{V.\,I.~Vernadsky Crimean Federal University, Simferopol, Russia. \email{fedyor@mail.ru}.}}
\date{}
\begin{document}
\begingroup
\makeatletter
\renewcommand{\@fnsymbol}[1]{\ensuremath{\ifcase#1\or *\or\spadesuit\or\text{\tiny\faMountain}\or\clubsuit\or\dagger\or\ddagger\fi}}
\let\origmakefntext\@makefntext
\long\def\@makefntext#1{\rightskip\z@\@plus3em\relax\origmakefntext{#1}}
\makeatother
\maketitle
\endgroup
\hypersetup{pdftitle={Descending Everest, Part III: Oracle Complexity of Monotone Variational Inequalities with (L0,L1)-Lipschitz Operators},pdfauthor={Mikhail Rudomanov, Ivan Fomin, Alexander V. Gasnikov, Fedor S. Stonyakin}}

\begin{abstract}
We consider monotone variational inequalities in which the Jacobian norm
of the operator $F\colon\R^d\to\R^d$ is bounded by an affine function
of the norm of the operator itself: $\norm{DF(z)}\le L_0+L_1\norm{F(z)}$.
The main case is the problem on the whole space, that is, finding a zero of $F$.
We assume that a solution exists and that an upper bound $R$ on the distance
from the starting point to the solution is known. Neither strong monotonicity
nor cocoercivity is assumed. Global Lipschitz continuity is not assumed either:
for $L_1>0$, the Jacobian norm may grow together with the norm of the operator.
We obtain upper and lower bounds on the number of oracle calls sufficient
to find a point $\hat z$ with $\norm{F(\hat z)}\le\varepsilon$. For $L_0R/\varepsilon\ge1$,
both bounds are, up to logarithmic factors in $L_1R$, of order
\[
 1+\frac{L_0R}{\varepsilon}+(L_1R)^{6/5}.
\]
The upper bound is attained by a deterministic extragradient-type method,
and the lower bound holds for arbitrary randomized methods. The same bounds
are obtained for the Minty gap on the ball of radius $R$, where
$L_0R/\varepsilon$ is replaced by $L_0R^2/\eta$.
In addition, we show that for a variational inequality on a bounded set
the numbers $L_0,L_1$ and the diameter of the set are not sufficient: a uniform
complexity bound requires additional information, for example, an upper bound
on the norm of the operator at the solution.
\end{abstract}

\begin{tcolorbox}[
  optimbox,
  colback=optimBlue!3!white,
  colframe=optimBlue!40!white,
  borderline west={1.4pt}{0pt}{optimBlue},
  top=7pt,bottom=7pt
]
\small
\noindent\textbf{Use of language models.}\enspace
The models \textit{Codex Astra 6} and \textit{Claude Opus 5.5} were used
in obtaining the theoretical results of this work.
They were also used to prepare and typeset the manuscript and to translate
the original text into English.
Nevertheless, every line of the paper, including every statement, every
proof, and every computation, has been personally edited and checked
by the authors. The authors bear sole responsibility for the correctness
of all mathematical statements, proofs, and conclusions.
\end{tcolorbox}

\section{Introduction}\label{sec:intro}

We consider the problem of finding a point $z\in\R^d$ such that $F(z)=0$, where $F\colon\R^d\to\R^d$ is a continuous monotone operator, that is, $\ip{F(u)-F(v)}{u-v}\ge0$ for all $u,v\in\R^d$.
A point $z^\ast\in Z$
is called a solution of the variational inequality on a convex set $Z$ if
$\ip{F(z^\ast)}{z-z^\ast}\ge0$ for all $z\in Z$. For $Z=\R^d$, the substitution
$z=z^\ast-tF(z^\ast)$, $t>0$, shows that this condition is equivalent to
the equality $F(z^\ast)=0$. Therefore, finding a zero of the operator is a variational
inequality on the whole space. The general theory of finite-dimensional variational inequalities
is presented in the monograph \citet{facchinei2003finite}.
On the whole space, a solution may fail to exist. The constant operator
$F\equiv c\ne0$ is monotone and satisfies any smoothness conditions, but it has
no zeros. On a compact set $Z$, a solution with the same operator exists: it is a
minimizer of the linear function $\ip cz$ on $Z$, that is, a solution of a linear
programming problem. Therefore, all bounds in this paper concern solvable
problems: we assume that a solution exists.

The main source of monotone operators is saddle-point problems. If a function
$\varphi(x,y)$ is convex in $x$ and concave in $y$, then the operator
$F(x,y)=\bigl(\nabla_x\varphi(x,y),-\nabla_y\varphi(x,y)\bigr)$ is monotone,
and its zeros are the saddle points of $\varphi$. Nash equilibrium conditions
in games where the loss function of each player is convex in the player's own
strategy are written in the same way: an equilibrium is a zero of the operator composed of the gradients of the players' loss
functions with respect to their strategies. If this operator is monotone, the game is called
monotone. A condition of this type (diagonal strict concavity) was introduced by
\citet{rosen1965existence}. Some problems of training generative adversarial
networks also reduce to monotone games \citep{gidel2019variational}. Another classical source is
traffic equilibria: the Wardrop equilibrium in the Beckmann model
is written as a variational inequality with the path cost operator,
which is monotone when the link cost functions are nondecreasing
\citep{dafermos1980traffic}. Methods for finding equilibria in the Beckmann
model and in the Nesterov--de~Palma stable dynamics model are considered
in \citet{kubentayeva2021finding}. The gradient of a convex function is a special
case of a monotone operator, but in general $F$ is not
a gradient. For $\varphi(x,y)=xy$, we obtain the rotation $F(x,y)=(y,-x)$,
and the step $z-hF(z)$ moves the point away from the solution for every $h>0$:
$\norm{z-hF(z)}^2=(1+h^2)\norm z^2$. Therefore, the main method for such problems
is the extragradient method \citep{korpelevich1976}, in which the operator
is evaluated at an extrapolation point, and its variant with one operator evaluation
per iteration \citep{popov1980}.

For a monotone operator with Lipschitz constant $L$, the complexity is well
studied. On a compact set $Z$, the accuracy of an approximate solution $z\in Z$ is usually
measured by the quantity $\sup_{u\in Z}\ip{F(u)}{z-u}$,
which is also called the Minty gap, or the restricted gap function
\citep{stonyakin2022generalized}. For a continuous monotone operator,
it is nonnegative and vanishes exactly at the solutions of the variational
inequality \citep{minty1962}.
In the Russian-language literature, the inequality $\sup_{u\in Z}\ip{F(u)}{z-u}\le\eta$
is used as an accuracy criterion for weak (Minty) variational inequalities
\citep{gasnikov2019adaptive}.
If the diameter of $Z$ equals $D$, the method \method{Mirror-Prox} finds a point with Minty
gap at most $\eta$ in $O(LD^2/\eta)$ iterations \citep{nemirovski2004},
and for linear operators this order cannot be improved \citep{nemirovsky1992}.
For the unconstrained problem, that is, the variational inequality on the whole
space, the natural accuracy criterion is the norm of the operator.
The Halpern iteration \citep{halpern1967} constructs a new point as a convex
combination of the starting point and a step of the method:
$z_{k+1}=\beta_kz_0+(1-\beta_k)T(z_k)$, where $T$ is a nonexpansive
map, for example, the resolvent of the operator $F$, and the weights $\beta_k$ decrease
to zero. The starting point keeps the iterates near itself, as an anchor holds
a ship. This is why it is called an anchor,
and such methods are called anchored methods. The meaning of the anchor is seen from the regularized
problem $F(z)+\beta(z-z_0)=0$, $\beta>0$. By the monotonicity of $F$, its solution
$z_\beta$ satisfies $\norm{z_\beta-z_0}\le\norm{z^\ast-z_0}$, hence
$\norm{F(z_\beta)}=\beta\norm{z_\beta-z_0}\le\beta\norm{z^\ast-z_0}$.
Anchored methods track $z_\beta$ with $\beta$ of order $L/k$, and the norm
of the operator at the last point decreases as $1/k$. The Halpern iteration made it possible to obtain
a bound of order $LR/\varepsilon$ up to a logarithm
\citep{diakonikolas2020}, and the accelerated anchored extragradient method
\method{EAG} (Extra Anchored Gradient) \citep{yoon2021} gives
$\norm{F(z_k)}^2=O(L^2R^2/k^2)$, that is, $O(LR/\varepsilon)$ oracle
calls. Here $R$ is the distance from the starting point to the solution.
For monotone operators, the optimal complexity scales linearly with $1/\varepsilon$,
whereas for convex minimization it is of order $\sqrt{LR^2/\varepsilon}$.

All these bounds rely on a global Lipschitz constant $L$,
which does not exist in many problems. For example, the function $e^x$ has an unbounded second
derivative, but it does not exceed the first one. Such growth, when the local
Lipschitz constant is bounded in terms of the norm of the operator itself, is described by
the $(L_0,L_1)$-smoothness condition. It was proposed for functions \citep{zhang2019gradient}
to explain the success of methods with gradient clipping: the local smoothness
constant can grow together with the gradient norm. In the first two parts of the
Descending Everest series of papers,
we found the oracle complexity of minimizing a single convex
$(L_0,L_1)$-smooth function \citep{rudomanov2026everest} and a sum of such
functions \citep{rudomanov2026everest2}. For operators, the analogous condition
$\norm{DF(z)}\le L_0+L_1\norm{F(z)}$
was considered in \citet{vankov2024}, where convergence of adaptive
projection and extragradient methods was proved under the additional
$p$-quasi-sharpness condition $\ip{F(u)}{u-u^\ast}\ge\mu\dist(u,U^\ast)^p$.
For a monotone $(L_0,L_1)$-Lipschitz operator without additional conditions,
\citet[Theorem~3.6]{choudhury2025} showed that the extragradient method
with stepsize inversely proportional to $L_0+L_1\norm{F(x_k)}$ finds a point
with $\norm{F}\le\varepsilon$ in $O\bigl(1+(L_0R/\varepsilon+L_1R)^2\bigr)$
oracle calls. The square arises because without an anchor the norm of the operator decreases as $1/\sqrt k$, and while the norm is large, the steps have length of order $1/L_1$. Subsequent developments include the extragradient method
with a moving anchor \citep{alcala2023moving}, \method{MAEG} with a local
stepsize test \citep{sun2026}, the method \method{PolyakEG} with Polyak-type
stepsizes \citep{yoon2026polyak}, the method \method{PF-AEG}, which does not require knowledge
of the parameters \citep{chu2026}, and an analysis of the extragradient method
and Popov's method in the presence of constraints \citep{greenberg2026}.
For a monotone operator without additional conditions, all these bounds
either have the same quadratic order or are not expressed in terms of
$L_0,L_1,R,\varepsilon$ alone. In the LongTails repository \citep{tyurin2026longtails}, a Lean-verified
bound $O(1+L_0R/\varepsilon+(L_1R)^{4/3})$ was obtained for the same class
(theorem \texttt{RadiusGeneralizedMonotone}).
We worked on this problem in parallel and independently and obtained a stronger bound
with exponent
$6/5$ on $L_1R$, which is optimal up to logarithmic factors:
it is matched by a lower bound
for arbitrary randomized methods. A summary table of known
bounds is given in Appendix~\ref{app:literature}.

The main difficulty in this setting is
that the distance to the solution gives only an exponential bound on the
norm of the operator at the starting point. Grönwall's lemma implies only the bound
$\norm{F(z_0)}\le(L_0/L_1)(e^{L_1R}-1)$, and it is attained. The length of a safe
step is inversely proportional to $L_0+L_1\norm{F(z)}$, so the method has to
pass through a region where the operator is exponentially large. For convex functions,
in the first part of the series, the cost of passing through this region turned out to be
$L_1R$ up to logarithms. Thus, continuing our series of papers,
we pose a natural question: how many oracle calls does it cost to pass through
this region for a monotone operator if only an upper
bound $R$ on the distance to the solution is known? In this paper, we answer it.

\paragraph{Notation.}
Throughout, $\norm{\cdot}$ denotes the Euclidean norm of a vector and the spectral norm
of a matrix, $[x,y]$ is the segment with endpoints $x,y\in\R^d$, $B(x,r)$ is the closed ball
of radius $r$ centered at $x$, $e_i$ is the $i$-th standard basis vector,
$\one$ is the all-ones vector, and $DF$ is the Jacobian matrix of the operator $F$.
The sum of the first $k$ terms of a sequence is denoted by the corresponding
capital letter, for example, $\Lambda_k=\lambda_1+\dots+\lambda_k$ and $\Lambda_0=0$.
The logarithm is natural, and $\log_+t=\max\{0,\log t\}$.

\paragraph{Main result.}
We show that for $L_0R/\varepsilon\ge1$ the minimax number of oracle calls
$\mathcal T^\ast$ required to find a point $\hat z$
with $\norm{F(\hat z)}\le\varepsilon$ satisfies the inequalities
\[
 c\Bigl(1+\frac{L_0R}{\varepsilon}+\frac{(L_1R)^{6/5}}{\log^{4/5}(2+L_1R)}\Bigr)
 \le1+\mathcal T^\ast
 \le C\Bigl(1+\frac{L_0R}{\varepsilon}+(L_1R)^{6/5}\Bigr)
\]
with absolute constants $c,C>0$ (Theorem~\ref{thm:main}).
The upper bound is attained by a deterministic method, and the lower bound holds
for arbitrary randomized methods, which may query the operator
at any points and output any point. The term $L_0R/\varepsilon$
is known from the Lipschitz case. The term $(L_1R)^{6/5}$ appears as the cost
of passing through the region of large operator values. Thus, for monotone operators,
this cost is larger than for gradients of convex functions, where it equals $L_1R$.
The same bounds with $L_0R/\varepsilon$ replaced by $L_0R^2/\eta$ are obtained
for the Minty gap computed over the ball of radius $R$ around the starting point.

\paragraph{Upper bound.}
The method \method{ABSEIL} consists of two stages. At the first stage (the method \method{GLISSADE}), while
$\norm{F(z)}>2L_0R\bigl(e^{(L_1R)^{2/5}}-1\bigr)/(L_1R)^{2/5}$, the method takes
steps of length $R/\bigl(2(L_1R)^{3/5}\bigr)$ along the direction
$-F(z)/\norm{F(z)}$. Near the solution, the operator is small by Grönwall's lemma, so
monotonicity guarantees that such a step decreases the squared distance to the solution
by $R^2/\bigl(4(L_1R)^{6/5}\bigr)$ (Lemma~\ref{lem:long-step}), and the first stage has
at most $4(L_1R)^{6/5}$ steps. At the second stage (the method \method{BELAY}), the extragradient
method with moving anchor \citep{alcala2023moving,sun2026} runs, in which the steps
are computed by an explicit formula. Its cost is of order
$L_0R/\varepsilon+(L_1R)^{6/5}$: the logarithm of $\norm{F}$ after the first stage
enters the bound under the square root (Proposition~\ref{prop:adaptive}). The step length
of the first stage is chosen so that the costs of the two stages have the same
order (Section~\ref{sec:total}).

\paragraph{Lower bound.}
The term $L_0R/\varepsilon$ is obtained by a reduction to solving linear systems, for which a lower bound for randomized methods was proved in \citet{derezinski2026}. The nonlinear term is obtained on a chain \citep{nemirovski1983problem,nesterov2018lectures}: each oracle answer
reveals at most one new coordinate. The worst-case operator is defined
as the inverse of $v\mapsto\partial(Q+I_K)(v)+Bv$, where $Q$ is a convex function on the weighted simplex $K$, and $B$ is a skew-symmetric block matrix. The generalized strong monotonicity of $\partial Q$ yields the $(L_0,L_1)$-Lipschitz continuity of the inverse operator, and the shift and the matrix $B$ force the coordinates to be revealed one at a time.
To prevent the method from guessing the unrevealed coordinates, the chain is randomly embedded into a high-dimensional space, and outside a large ball the operator is replaced by a known linear one
\citep{woodworth2016finite,carmon2020lower}.

\paragraph{Bounded set.}
For a variational inequality on a compact set, a solution no longer has to be a zero
of the operator. For each continuous operator, the norm of $F(z^\ast)$ is finite,
but it is not uniformly bounded over the class of operators with given $L_0$, $L_1$,
and diameter of the set.
We show that, as a consequence, the minimax complexity of this class is infinite
(Theorem~\ref{thm:constrained-impossible}),
and a uniform bound requires an additional parameter, for example,
an upper bound on $\norm{F(z^\ast)}$.

\paragraph{Organization of the paper.}
The problem setting and the main result are given in Section~\ref{sec:setting},
the method and the upper bound in Section~\ref{sec:upper}, and the lower bounds
in Section~\ref{sec:lower}. Section~\ref{sec:constraints} is devoted to
variational inequalities on a bounded set, and
Section~\ref{sec:conclusion} to open questions. The proofs are deferred
to Appendices~\ref{app:tools}--\ref{app:constraints}, and Appendix~\ref{app:literature}
contains a summary table of known results.

\section{Problem setting and main result}\label{sec:setting}

We consider an operator $F\colon\R^d\to\R^d$. We are given numbers $L_0,R>0$,
$L_1\ge0$, and a starting point $z_0\in\R^d$. All general assumptions are collected below.

\begin{assumption}[monotonicity]\label{asm:monotone}
The operator $F$ is continuous and $\ip{F(u)-F(v)}{u-v}\ge0$ for all $u,v\in\R^d$.
\end{assumption}

\begin{assumption}[$(L_0,L_1)$-Lipschitz continuity]\label{asm:smooth}
For all $u,v\in\R^d$,
\begin{equation}\label{eq:smooth}
 \norm{F(u)-F(v)}\le\Bigl(L_0+L_1\sup_{w\in[u,v]}\norm{F(w)}\Bigr)\norm{u-v}.
\end{equation}
\end{assumption}

\begin{assumption}[solution]\label{asm:root}
There exists a point $z^\ast$ such that $F(z^\ast)=0$ and $\norm{z^\ast-z_0}\le R$.
\end{assumption}

Assumption~\ref{asm:root} singles out the class of solvable problems, and all bounds
of this paper, upper and lower, refer to this class. Without it, the problem may
have no solution, as for the constant operator $F\equiv c\ne0$.

Assumption~\ref{asm:smooth} is the symmetric $(L_0,L_1)$-smoothness condition
\citep[Assumption~1.3]{gorbunov2024methods}, written for an operator instead of
a gradient. For operators, it is called $(L_0,L_1)$-Lipschitz continuity
\citep{choudhury2025,greenberg2026}. It coincides with the $\alpha$-symmetric
generalized smoothness condition of \citet{vankov2024} for $\alpha=1$.
Inequality~\eqref{eq:smooth} is required on the whole space, since
the method may query the operator outside the ball $B(z_0,R)$. We assume
neither strong monotonicity of the operator nor that it is a potential operator,
that is, the gradient of some function. The following lemma shows
that for continuously differentiable operators, condition~\eqref{eq:smooth}
is equivalent to the Jacobian norm bound $\norm{DF(z)}\le L_0+L_1\norm{F(z)}$,
by which $(L_0,L_1)$-Lipschitz continuity of an operator is defined in the introduction.

\proofdown{lem:class}
\begin{lemma}[characterization of the class]\label{lem:class}
\leavevmode
\begin{enumerate}
\item[(a)] A continuous operator $F$ satisfies~\eqref{eq:smooth} if
and only if for any $x,v\in\R^d$ the map $t\mapsto F(x+tv)$
is locally absolutely continuous and for almost all $t$
\[
 \Bigl\|\frac{d}{dt}F(x+tv)\Bigr\|\le\bigl(L_0+L_1\norm{F(x+tv)}\bigr)\norm v.
\]
\item[(b)] If $F\in C^1(\R^d)$, then $F$ satisfies~\eqref{eq:smooth} if
and only if $\norm{DF(z)}\le L_0+L_1\norm{F(z)}$ for all $z\in\R^d$.
\item[(c)] If $F$ is locally Lipschitz and $\norm{DF(z)}\le L_0+L_1\norm{F(z)}$
for almost all $z$, then $F$ satisfies~\eqref{eq:smooth}.
\end{enumerate}
\end{lemma}

The main consequence of condition~\eqref{eq:smooth} is a bound on the rate of change
of the operator along a segment. It replaces Lipschitz continuity in all proofs
of this paper.

\proofdown{lem:gronwall}
\begin{lemma}[consequences of Grönwall's lemma]\label{lem:gronwall}
Let $F$ satisfy~\eqref{eq:smooth}. Then for any $u,v\in\R^d$
\begin{align*}
 \text{(a)}\quad L_0+L_1\norm{F(v)}&\le\bigl(L_0+L_1\norm{F(u)}\bigr)e^{L_1\norm{u-v}},\\
 \text{(b)}\quad \norm{F(v)-F(u)}&\le\bigl(L_0+L_1\norm{F(u)}\bigr)\frac{e^{L_1\norm{u-v}}-1}{L_1}\\
 &\le\bigl(L_0+L_1\norm{F(u)}\bigr)\norm{u-v}e^{L_1\norm{u-v}}.
\end{align*}
\end{lemma}

For $L_1=0$, the fraction $(e^{L_1\norm{u-v}}-1)/L_1$ in Lemma~\ref{lem:gronwall}
is understood as $\norm{u-v}$. The proof of the lemma uses condition~\eqref{eq:smooth}
only on the segment $[u,v]$. If the map
$t\mapsto F(u+t(v-u))$ is absolutely continuous on $[0,1]$ and the bound
of Lemma~\ref{lem:class}(a) with some constant $A>0$ in place of $L_0$ holds
only for points of the segment $[u,v]$, then (a) and (b) hold with $A$ in place of $L_0$.

In particular, by Assumption~\ref{asm:root}, Lemma~\ref{lem:gronwall}(b)
with $u=z^\ast$ gives $\norm{F(z_0)}\le(L_0/L_1)(e^{L_1R}-1)$. This bound
is attained. For the convex function
$f(z)=L_0L_1^{-2}\bigl(e^{L_1\norm z}-1-L_1\norm z\bigr)$, the operator
$F=\nabla f$ satisfies Assumptions~\ref{asm:monotone}--\ref{asm:root}
with $z^\ast=0$ \citep{rudomanov2026everest}, and
$\norm{F(z)}=(L_0/L_1)(e^{L_1\norm z}-1)$.
Thus, the norm of the operator at the starting point may be exponentially
large in $L_1R$.

\paragraph{Oracle model.}
The algorithm knows $L_0$, $L_1$, $R$, $z_0$, and the required accuracy. It obtains
information about the operator only from an oracle that, for a query $z\in\R^d$,
returns the exact value $F(z)$. Queries are allowed anywhere in the space.
The query points, the stopping time $\tau$, and the output point $\hat z$ are
measurable functions of the received oracle answers and the random variables
of the algorithm. Complexity is measured by the number of oracle calls
\citep[Section~1.1]{nesterov2018lectures}. The word ``randomized'' refers
to the algorithm: the oracle answers contain no noise.

\paragraph{Accuracy criteria.}
The main criterion is the norm of the operator at the output point: $\norm{F(\hat z)}\le\varepsilon$.
For comparison with the theory of variational inequalities, we also consider the Minty gap
\citep{minty1962} computed over the ball $B(z_0,R)$:
\begin{equation}\label{eq:gap}
 \Gap_R(z)=\sup_{u\in B(z_0,R)}\ip{F(u)}{z-u}.
\end{equation}
Here the ball only defines the set over which the supremum is taken: the queries
and the output point may lie outside it. Since $z^\ast\in B(z_0,R)$
and $\ip{F(z^\ast)}{z-z^\ast}=0$, the gap is nonnegative. The supremum cannot be taken
over the whole space: for $F(z)=Sz$ with a nonzero skew-symmetric
matrix $S$, we have $\sup_{u\in\R^d}\ip{Su}{z-u}=\sup_{u\in\R^d}\ip{Su}z=+\infty$
if $Sz\ne0$. A small norm of the operator implies a small $\Gap_R(z)$.

\proofdown{lem:gap-norm}
\begin{lemma}[gap and the norm of $F$]\label{lem:gap-norm}
Under Assumptions~\ref{asm:monotone} and~\ref{asm:root}, for any point $z\in\R^d$
we have $0\le\Gap_R(z)\le\bigl(\norm{z-z_0}+R\bigr)\norm{F(z)}$.
\end{lemma}

\paragraph{Oracle complexity.}
We denote by $\mathcal P_d$ the set of operators $F\colon\R^d\to\R^d$
satisfying Assumptions~\ref{asm:monotone}--\ref{asm:root}.
Let $\tau_{\mathsf M}(F)$ be the number of oracle calls made by the algorithm $\mathsf M$
on the operator $F$, and let $\hat z_{\mathsf M}(F)$ be its output point.

\begin{definition}[oracle complexity]\label{def:complexity}
An algorithm is called admissible if for any $d$ and $F\in\mathcal P_d$
it stops almost surely and
\[
 \Pr\bigl\{\norm{F(\hat z_{\mathsf M}(F))}\le\varepsilon\bigr\}\ge\frac23.
\]
The complexity of an algorithm and the minimax complexity are defined as
\begin{equation}\label{eq:minimax}
\begin{gathered}
 \mathcal T(\mathsf M)=\sup_{d\ge1}\sup_{F\in\mathcal P_d}\E\tau_{\mathsf M}(F),
 \\
 \mathcal T^\ast(\varepsilon)=\inf_{\mathsf M\ \text{admissible}}\mathcal T(\mathsf M).
\end{gathered}
\end{equation}
\end{definition}

We denote by $\mathcal T^\ast_{\det}(\varepsilon)$ the same quantity where
the infimum is taken over deterministic algorithms for which
$\norm{F(\hat z_{\mathsf M}(F))}\le\varepsilon$ on every operator.
The quantities $\mathcal T^\ast_{\Gap}(\eta)$ and $\mathcal T^\ast_{\Gap,\det}(\eta)$
are defined in the same way with the condition $\norm{F(\hat z)}\le\varepsilon$
replaced by $\Gap_R(\hat z)\le\eta$.

Clearly, $\mathcal T^\ast\le\mathcal T^\ast_{\det}$. The class $\mathcal P_d$
is invariant under shifts, so in the lower bounds we assume $z_0=0$.

\begin{theorem}[two-sided complexity bound]\label{thm:main}
There exist absolute constants $c,C>0$ such that for
$L_0R/\varepsilon\ge1$
\begin{equation}\label{eq:main}
 c\Bigl(1+\frac{L_0R}{\varepsilon}+\frac{(L_1R)^{6/5}}{\log^{4/5}(2+L_1R)}\Bigr)
 \le1+\mathcal T^\ast(\varepsilon)\le1+\mathcal T^\ast_{\det}(\varepsilon)
 \le C\Bigl(1+\frac{L_0R}{\varepsilon}+(L_1R)^{6/5}\Bigr),
\end{equation}
and for $L_0R^2/\eta\ge1$
\begin{equation}\label{eq:main-gap}
 c\Bigl(1+\frac{L_0R^2}{\eta}+\frac{(L_1R)^{6/5}}{\log^{4/5}(2+L_1R)}\Bigr)
 \le1+\mathcal T^\ast_{\Gap}(\eta)\le1+\mathcal T^\ast_{\Gap,\det}(\eta)
 \le C\Bigl(1+\frac{L_0R^2}{\eta}+(L_1R)^{6/5}\Bigr).
\end{equation}
\end{theorem}

The lower bounds remain valid on the subclass of operators from $C^2(\R^d)$.
The upper bounds are proved in Theorem~\ref{thm:upper}. They are attained by the method \method{ABSEIL}
(Section~\ref{sec:method}), which outputs a point from the ball $B(z_0,5R)$
and uses the known $L_0$, $L_1$, $R$. The term $L_0R/\varepsilon$ (respectively, $L_0R^2/\eta$)
in the lower bounds is obtained in Proposition~\ref{prop:linear-lower}, and the nonlinear
term in Theorem~\ref{thm:nonlinear-lower}. The two lower bounds
are combined by the inequality $\max\{x,y\}\ge(x+y)/2$.

The upper and lower bounds differ only by the factor $\log^{4/5}(2+L_1R)$
in the nonlinear term. In particular, if $L_0R/\varepsilon\ge(L_1R)^{6/5}$,
then $1+\mathcal T^\ast$ and $1+\mathcal T^\ast_{\det}$ are equal to $\Theta(1+L_0R/\varepsilon)$,
and for $L_1=0$ we recover the classical bound for Lipschitz operators.
The condition $L_0R/\varepsilon\ge1$ is essential for the lower bound. If, for example,
$(L_0/L_1)(e^{L_1R}-1)\le\varepsilon$, then by Lemma~\ref{lem:gronwall} the point
$z_0$ is acceptable without any oracle calls.

All three quantities $L_0R/\varepsilon$, $L_0R^2/\eta$, and $L_1R$ are dimensionless.
If $\widetilde F(x)=\gamma F(z_0+\beta x)$ with $\gamma,\beta>0$, then
$\widetilde F$ satisfies the assumptions with the parameters
$\gamma\beta L_0$, $\beta L_1$, $R/\beta$, and the accuracies become
$\gamma\varepsilon$ and $\gamma\eta/\beta$. Therefore, the complexity measured by $\norm{F}$
can depend only on $L_0R/\varepsilon$ and $L_1R$, and the complexity
measured by the gap only on $L_0R^2/\eta$ and $L_1R$.

\section{Upper bound}\label{sec:upper}

In this section, we construct a method that attains the upper bound
of Theorem~\ref{thm:main}. Throughout the section,
Assumptions~\ref{asm:monotone}--\ref{asm:root} hold.

\subsection{The method}\label{sec:method}

The method \method{ABSEIL} consists of two stages\footnote{The names of the methods are taken
from mountaineering: an abseil is a descent on a rope, a glissade is a slide down a snow
slope, and a belay is securing a climber with a rope.}. At the first stage, while the norm of the operator is large,
the method \method{GLISSADE} (Algorithm~\ref{alg:glissade}) takes long
steps along the direction $-F(z)/\norm{F(z)}$. At the second stage, the method
\method{BELAY} (Algorithm~\ref{alg:belay}) runs. It is an extragradient method with a moving
anchor. The moving anchor was proposed by \citet{alcala2023moving},
and for $(L_0,L_1)$-Lipschitz operators such a method, called \method{MAEG},
was studied by \citet{sun2026}. In \method{MAEG}, a trial step is accepted after
checking the stepsize condition at already computed points and is decreased if the condition
fails, so rejected steps require additional oracle
calls. In \method{BELAY}, the stepsize is computed before the oracle call from already known
values of the operator, and this stepsize always satisfies the required condition
(Proposition~\ref{prop:adaptive}).

The method \method{ABSEIL} works as follows. If $L_1R\ge1$, it runs \method{GLISSADE} from the point
$z_0$ and obtains a point $s$. If $L_1R<1$, it sets $s=z_0$ and queries
$F(z_0)$. If $\norm{F(s)}\le\varepsilon$, the method outputs $s$. Otherwise, it runs
\method{BELAY} from the point $s$ and outputs the resulting point.

\begin{algorithm}[ht]
\caption{\method{GLISSADE} (first stage)}\label{alg:glissade}
\begin{algorithmic}[1]
\Require $z_0$, $L_0$, $L_1$, $R$, $\varepsilon$, where $L_1R\ge1$
\State $z\gets z_0$
\State query $F(z)$
\While{$\norm{F(z)}>\max\bigl\{\varepsilon,\ 2L_0R\bigl(e^{(L_1R)^{2/5}}-1\bigr)/(L_1R)^{2/5}\bigr\}$}
 \State $z\gets z-\dfrac{R}{2(L_1R)^{3/5}}\,\dfrac{F(z)}{\norm{F(z)}}$
 \State query $F(z)$
\EndWhile
\State\Return $z$
\end{algorithmic}
\end{algorithm}

\begin{algorithm}[ht]
\caption{\method{BELAY} (second stage)}\label{alg:belay}
\begin{algorithmic}[1]
\Require $s$, $F(s)$, $L_0$, $L_1$, $\varepsilon$
\State $u_0\gets s$, $y_0\gets s$, $\Lambda_0\gets0$
\For{$k=1,2,\ldots$}
 \If{$k=1$}
  \State $\lambda_1\gets\dfrac1{6\bigl(L_0+L_1\norm{F(s)}\bigr)}$
 \Else
  \State $B_k\gets\norm{F(y_{k-1})}+\dfrac{3\norm{u_{k-1}-y_{k-1}}}{2\Lambda_{k-1}}$
  \State $\lambda_k\gets\min\Bigl\{4\lambda_{k-1},\ \dfrac1{6(L_0+L_1B_k)}\Bigr\}$
 \EndIf
 \State $\Lambda_k\gets\Lambda_{k-1}+\lambda_k$
 \State $\alpha_k\gets\dfrac{3\lambda_k}{2\Lambda_{k-1}+3\lambda_k}$
 \State $v_k\gets\alpha_ku_{k-1}+(1-\alpha_k)y_{k-1}$
 \State $x_k\gets v_k-(1-\alpha_k)\lambda_kF(y_{k-1})$
 \State query $F(x_k)$ \Comment{for $k=1$ we have $x_1=s$, so no query is needed}
 \State $y_k\gets v_k-\lambda_kF(x_k)$
 \State query $F(y_k)$
 \State $u_k\gets u_{k-1}-\tfrac16\lambda_kF(y_k)$
 \If{$\norm{F(y_k)}\le\varepsilon$}
  \State\Return $y_k$
 \EndIf
\EndFor
\end{algorithmic}
\end{algorithm}

\subsection{Long steps}\label{sec:long-steps}

Long steps are based on the following observation. Let $\norm{x-z^\ast}\le R$
and let the norm $\norm{F(x)}$ be large. Consider the point
$w=z^\ast+R(L_1R)^{-3/5}F(x)/\norm{F(x)}$, shifted from the solution
in the direction $F(x)$. It lies at a distance $R/(L_1R)^{3/5}$
from the solution, so by Lemma~\ref{lem:gronwall}(b) with $u=z^\ast$
\[
 \norm{F(w)}\le\frac{L_0}{L_1}\bigl(e^{(L_1R)^{2/5}}-1\bigr),
\]
and this quantity is small compared to $\norm{F(x)}$. By monotonicity,
$\ip{F(x)}{x-w}\ge\ip{F(w)}{x-w}\ge-\norm{F(w)}\norm{x-w}$.
Hence the projection of the vector $x-z^\ast$ onto the direction $F(x)/\norm{F(x)}$
is at least $R/(L_1R)^{3/5}$ up to a small correction. Therefore, a step
of length $R/\bigl(2(L_1R)^{3/5}\bigr)$ along $-F(x)/\norm{F(x)}$ brings the point closer
to the solution, although the solution itself is unknown to the method. The precise statement
is given in the following lemma.

\proofdown{lem:long-step}
\begin{lemma}[long step]\label{lem:long-step}
Let $\norm{x-z^\ast}\le R$ and $L_1R\ge1$. If
$\norm{F(x)}>2L_0R\bigl(e^{(L_1R)^{2/5}}-1\bigr)/(L_1R)^{2/5}$, then
\[
 \Bigl\|x-\frac{R}{2(L_1R)^{3/5}}\frac{F(x)}{\norm{F(x)}}-z^\ast\Bigr\|^2
 \le\norm{x-z^\ast}^2-\frac{R^2}{4(L_1R)^{6/5}}.
\]
\end{lemma}

Consequently, \method{GLISSADE} makes at most $4(L_1R)^{6/5}$ steps, all its points
remain in the ball $B(z^\ast,R)$, and the point $s$ at which the second stage starts
satisfies the inequalities
\begin{equation}\label{eq:stage-one}
\begin{gathered}
 \norm{s-z^\ast}\le R,\\ \norm{s-z_0}\le2R,\\
 \norm{F(s)}\le2L_0R\bigl(e^{(L_1R)^{2/5}}-1\bigr)/(L_1R)^{2/5}.
\end{gathered}
\end{equation}
If $L_1R<1$, the first stage is skipped and $s=z_0$.

\subsection{Second stage: the method \texorpdfstring{\method{BELAY}}{BELAY}}\label{sec:maeg}

In the method \method{BELAY}, at each iteration the anchor $u_k$
is shifted by $\lambda_kF(y_k)/6$ in the direction $-F(y_k)$, and the extragradient
step is made from the point $v_k$, which is a convex combination of the anchor $u_{k-1}$
and the previous point $y_{k-1}$ with the weight $\alpha_k$. The analysis relies
on two cumulative bounds. The first one is the bound for the Lyapunov function
\[
 \ip{u_k-y_k}{F(y_k)}-\frac{\Lambda_k}3\norm{F(y_k)}^2
\]
from \citet[Proposition~2, Theorem~1]{sun2026}. It controls the norm
of the operator via the distance from the anchor to the solution. The second bound controls
the distance $\norm{u_k-y_k}$ from the current point to the anchor. It is this bound that allows
us to choose the stepsize from the computable quantity $B_k$, which bounds the length of the trial step:
$\norm{x_k-y_{k-1}}\le\lambda_kB_k$.

\proofdown{lem:energies}
\begin{lemma}[two cumulative bounds]\label{lem:energies}
Let the points of Algorithm~\ref{alg:belay} be constructed with
arbitrary $\lambda_k>0$ such that
$\lambda_k\le4\lambda_{k-1}$ for $k\ge2$, and let
$\norm{s-z^\ast}\le R$. If for all $i\le k$
\begin{equation}\label{eq:accept}
 \lambda_i\norm{F(x_i)-F(y_i)}\le\frac14\norm{x_i-y_i},
\end{equation}
then
\begin{gather*}
 \norm{u_k-z^\ast}\le R,\\
 \norm{F(y_k)}\le\frac{3R}{\Lambda_k},\\
 \sum_{i=1}^k\lambda_i\Lambda_i\norm{F(y_i)}^2\le9R^2,\\
 \norm{u_k-y_k}^2+5\norm{u_k-z^\ast}^2+\sum_{i=1}^k\alpha_i\norm{u_{i-1}-y_{i-1}}^2
 \le\frac{37}3R^2.
\end{gather*}
\end{lemma}

Condition~\eqref{eq:accept} is the usual stepsize condition of the extragradient
method. The following proposition shows that the stepsizes of \method{BELAY}
satisfy it and bounds the number of its iterations. The ratio
$\lambda_k/\lambda_1$ is not bounded here: the stepsize may both increase
and decrease, but by at most a factor of four per iteration.

\proofdown{prop:adaptive}
\begin{proposition}[number of iterations of the second stage]\label{prop:adaptive}
Let $\norm{s-z^\ast}\le R$ and
\[
 \mathcal H=\log\frac{\bigl(18L_0R/\varepsilon+1\bigr)\bigl(L_0+L_1\norm{F(s)}\bigr)}{L_0}.
\]
Then the stepsizes of Algorithm~\ref{alg:belay} started from the point $s$ satisfy
condition~\eqref{eq:accept}, and the algorithm stops after $K$ iterations, where
\[
 K\le3+\frac{36L_0R}{\varepsilon}+330L_1R\sqrt{\mathcal H}+\mathcal H.
\]
The output point satisfies $\norm{F(y_K)}\le\varepsilon$
and $\norm{y_K-z^\ast}<4R$.
\end{proposition}

The quantity $\mathcal H$ is the logarithm of the ratio of the largest stepsize that
may be needed to the initial stepsize $\lambda_1$. Without the first stage, it may
be of order $L_1R$, and then the term $L_1R\sqrt{\mathcal H}$ is of order
$(L_1R)^{3/2}$. The first stage reduces $\norm{F(s)}$ to $2L_0R\bigl(e^{(L_1R)^{2/5}}-1\bigr)/(L_1R)^{2/5}$,
and then $\mathcal H\lesssim1+(L_1R)^{2/5}+\log(2+L_0R/\varepsilon)+\log(2+L_1R)$.

\subsection{The final bound}\label{sec:total}

Let us explain the choice of the step length at the first stage. Lemma~\ref{lem:long-step}
and its proof remain valid if $(L_1R)^{2/5}$ is replaced by any
number $b\in[1,L_1R]$: a step of length $b/(2L_1)$ made when
$\norm{F(x)}>2L_0R(e^b-1)/b$ reduces the squared distance to the solution
by $b^2/(4L_1^2)$. Then the first stage costs at most $4(L_1R)^2/b^2$ oracle
calls. After it, $\mathcal H\lesssim1+b+\log(2+L_0R/\varepsilon)+\log(2+L_1R)$,
and by Proposition~\ref{prop:adaptive} the second stage costs of order
$L_0R/\varepsilon+L_1R\sqrt b$ up to lower-order terms. The first
quantity decreases in $b$, the second one increases, and for $b=(L_1R)^{2/5}$
both are equal to $(L_1R)^{6/5}$.

\proofdown{thm:upper}
\begin{theorem}[upper bound]\label{thm:upper}
The method \method{ABSEIL} outputs a point $\hat z$ such that
$\norm{F(\hat z)}\le\varepsilon$ and $\norm{\hat z-z_0}<5R$ after
\[
 O\Bigl(1+\frac{L_0R}{\varepsilon}+(L_1R)^{6/5}\Bigr)
\]
oracle calls. If it is run with the accuracy $\varepsilon=\eta/(6R)$,
then $\Gap_R(\hat z)\le\eta$, and the number of oracle calls is
$O\bigl(1+L_0R^2/\eta+(L_1R)^{6/5}\bigr)$.
\end{theorem}

The method uses no projections, no line search, and no auxiliary subproblems.
Every computed stepsize is accepted, and all oracle calls, including
the evaluation of the operator at the output point, are counted in the bound. Theorem~\ref{thm:nonlinear-lower}
shows that the exponent $6/5$ cannot be decreased.

\section{Lower bounds}
\label{sec:lower}
In the model of Definition~\ref{def:complexity}, the algorithm queries values $F(u)$ and must find a point where $F$ is almost zero. For the lower bound, we construct an operator whose answers for a long time do not give enough information to find such a point with probability at least $2/3$.

For $L_0R/\varepsilon\ge1$, the complexity in the paper has two different sources:
\[
1+\mathcal T^\ast(\varepsilon)\ \gtrsim\
1+\frac{L_0R}{\varepsilon}
+\frac{(L_1R)^{6/5}}{\log^{4/5}(2+L_1R)}.
\]
We first construct a chain and explain the origin of the nonlinear term. We then transfer the bound to the Minty gap, obtain the linear term, and consider a bounded set.

All constructed operators have a zero inside the ball of radius $R$ (Section~\ref{sec:lower-radius}), that is, the problems are solvable in the sense of Assumption~\ref{asm:root}. The difficulty in the lower bound comes not from the absence of a solution but from the fact that it cannot be found for a long time.

The worst-case operator must satisfy three requirements.
\begin{itemize}
\item The operator has a solution $z^*$, $F(z^*)=0$, with $\norm{z^*}\le R$.
\item The operator is monotone: $\ip{F(u)-F(u')}{u-u'}\ge0$.
\item The operator is $(L_0,L_1)$-Lipschitz in the sense of~\eqref{eq:smooth}. By Lemma~\ref{lem:class}, for a smooth operator this is equivalent to $\norm{DF(u)}\le L_0+L_1\norm{F(u)}$, that is, the Jacobian norm may be large where the norm of the operator is large.
\end{itemize}

We distinguish two spaces: $u$ is a query of the algorithm, and $v=\Phi(u)$ is the oracle answer. We will restrict the geometry of the answers $v$, although the algorithm chooses the queries $u$ freely.

The worst-case operator is chosen before the random decisions of the algorithm. The success requirement is probability at least $2/3$ on every operator of the class. The complexity $\mathcal T^\ast$ is measured by the worst-case expected number of queries, and the dimension is unbounded.

Throughout this section, we move the starting point to the origin, that is, we set $z_0=0$.

\subsection{What the chain must hide}
\label{sec:lower-target}
We construct an operator whose answers reveal coordinates only one at a time: if the query does not yet use the coordinates after $k$, the answer can reveal at most coordinate $k+1$. This is a property of the constructed operator, not a restriction on the algorithm. The queries $u$ remain arbitrary, and against arbitrary randomized algorithms the property is transferred by the random embedding of Section~\ref{sec:hiding}.

Let $\theta>0$ be a smallness threshold. After choosing the parameters, we want to obtain
\begin{align*}
|u_i|\le\theta\ \text{for}\ i\ge j
&\ \Longrightarrow\ \Phi_i(u)=0\ \text{for}\ i>j,\\
|u_n|\le\theta
&\ \Longrightarrow\ \norm{\Phi(u)}>2\varepsilon.
\end{align*}
Moreover, the values of the nonzero coordinates of the answer must also depend only on $u_1,\ldots,u_j$.

Thus, the chain must have two properties. First, each answer reveals at most one new coordinate. Second, as long as the last coordinate is not revealed, the norm of the answer is greater than $2\varepsilon$. In Section~\ref{sec:chain}, we construct an operator with the second property, and in Section~\ref{sec:lower-face} we verify the first.

\subsection{Geometry of the answers and the inverse operator}\label{sec:chain}

\subsubsection{Weighted simplex}
\label{sec:lower-simplex}
Let $G_1\ge\cdots\ge G_n>0$. We force all answers to lie in the weighted simplex
\[
K=\left\{v\ge0:\sum_{i=1}^n\frac{v_i}{G_i}\le1\right\}.
\]
This is the convex hull of $0,G_1e_1,\ldots,G_ne_n$. Zero belongs to this simplex but not to its outer face $\sum_i v_i/G_i=1$.

\begin{lemma}[outer face]\label{lem:face}
If $v\ge0$ and $\sum_iv_i/G_i=1$, then $\norm v\ge G_n/\sqrt n$.
\end{lemma}
\noindent Indeed,
\[
1=\sum_i\frac{v_i}{G_i}\le\frac1{G_n}\sum_i v_i
\le\frac{\sqrt n}{G_n}\norm v.
\]

The task is now concrete: as long as $u_n$ is small, we must keep the answer $v$ on the outer face of the simplex.

\subsubsection{Inverse operator}
\label{sec:lower-inverse}
In the first part of the series~\citep{rudomanov2026everest}, the worst-case functions were constructed via conjugation. Here we use an analogous technique of inverting a monotone operator. It is hard to give an explicit formula $v=\Phi(u)$ with all the required properties. It is more convenient to describe, for each admissible answer $v\in K$, the queries $u$ for which the oracle returns this answer:
\begin{equation}\label{eq:inverse-relation}
u+\sigma\in\partial(Q+I_K)(v)+Bv.
\end{equation}
Here $Q$ is convex, $B^\top=-B$, $\sigma$ is a shift, and $I_K$ is the indicator function of the set $K$: $I_K=0$ on $K$ and $I_K=+\infty$ outside $K$.

We denote the $(L_0,L_1)$-Lipschitz constants of the intermediate operator by $\ell_0,\ell_1>0$. At the end, we choose them based on $L_0,L_1$.

At a smooth point, \eqref{eq:inverse-relation} means $u+\sigma=\nabla Q(v)+Bv+y$ with $y\in N_K(v)$. The normal cone $N_K(v)$ consists of the $y$ for which $\ip y{w-v}\le0$ for all $w\in K$: these are the outer normals to the feasible directions. In the interior of $K$, only $y=0$ remains, and additional normals appear on the boundary.

\subsubsection{The function \texorpdfstring{$Q$}{Q}}
\label{sec:lower-Q}
For $v\ge0$, $v\ne0$, set $s=\sum_i v_i$. Define $Q$ as
\begin{equation}\label{eq:Q}
Q(v)=\underbrace{\frac2{\ell_1}\int_0^s\log\left(1+\frac{\ell_1t}{\ell_0}\right)dt}_{\text{change of the total size}}
+\underbrace{\frac4{\ell_1}\norm v\log\frac{\sqrt n\norm v}{s}}_{\text{change of the distribution}},
\end{equation}
and $Q(0)=0$.

Decompose $\xi\in\R^n$ as $\xi=\alpha v/s+\xi_\perp$, where $\alpha=\one^\top\xi$ and $\one^\top\xi_\perp=0$. The computation in Appendix~\ref{app:lower-chain} gives
\[
D^2Q(v)[\xi,\xi]
\ge\frac{2\alpha^2}{\ell_0+\ell_1s}
+\frac{4\norm{\xi_\perp}^2}{\ell_1\norm v}
\ge\frac{\norm\xi^2}{\ell_0+\ell_1\norm v}.
\]
Integration along a segment turns this bound into generalized strong monotonicity of the subdifferential: for $a\in\partial(Q+I_K)(v)$ and $a'\in\partial(Q+I_K)(v')$,
\[
\ip{a-a'}{v-v'}\ge\frac{\norm{v-v'}^2}{\ell_0+\ell_1\max\{\norm v,\norm{v'}\}}.
\]
The skew-symmetric $B$ does not change it, since $\ip{Bh}h=0$. After inversion, it yields monotonicity and $(\ell_0,\ell_1)$-Lipschitz continuity of $\Phi$.

\begin{remark}
In one variable, for $u=Q'(v)$, we have
\[
\frac{dv}{du}=\frac1{Q''(v)}.
\]
Therefore, the lower bound on the curvature $Q''(v)\ge1/(\ell_0+\ell_1v)$ turns into the upper bound $\Phi'(u)\le\ell_0+\ell_1\Phi(u)$.
\end{remark}

The order of the coordinates is set by the matrix $B$ and the shift $\sigma$.

\subsubsection{Blocks}
\label{sec:lower-blocks}
Let $H\ge2$ and $m\ge1$ be integers. We split the $n=Hm$ coordinates into $H$ consecutive blocks of $m$. Denote by $\bl(i)$ the index of the block of coordinate $i$, and by $\pos(i)$ its index within the block. Within a block the weight is the same, and between blocks it decreases geometrically:
\[
G_i=\frac{\ell_0}{\ell_1}n^3e^{H-\bl(i)}.
\]
The large $G_n=(\ell_0/\ell_1)n^3$ ensures a large error even on the last block. The geometric decrease in the weights allows us to enforce the order in which the blocks are revealed.

The matrix $B$ is defined as follows:
\[
B_{ij}=
\begin{cases}
\dfrac{5\sqrt m\log(e\sqrt n)}{\ell_1G_i}\operatorname{sign}(i-j),&\bl(i)=\bl(j),\\[2ex]
0,&\bl(i)\ne\bl(j).
\end{cases}
\]
Since $G_i=G_j$ within a block, the matrix is skew-symmetric.

Below we choose a shift $\sigma_i$ that decreases within a block. The matrix $B$ links the coordinates within a block. The decrease of $\sigma_i$ gives an advantage to the earlier coordinate. The different weights $G_i$ forbid transitions between blocks. The blocks are needed so that the cost of the internal order depends on $\sqrt m$ rather than on $\sqrt n$.

\subsection{Why the chain is revealed one coordinate at a time}
\label{sec:lower-face}
For $v\in K$, $v\ne0$, differentiating $Q$ gives $\nabla_iQ(v)=\vartheta(v)+\gamma(v)v_i$, where
\begin{equation}\label{eq:perspective-gradient}
\begin{gathered}
\vartheta=\frac2{\ell_1}\log\left(1+\frac{\ell_1s}{\ell_0}\right)-\frac{4\norm v}{\ell_1s},
\\
\gamma=\frac4{\ell_1\norm v}\left(1+\log\frac{\sqrt n\norm v}{s}\right).
\end{gathered}
\end{equation}
The common term $\vartheta$ is the same for all coordinates. On $K$, we have $\norm v\le s\le G_1$ and $1\le\sqrt n\norm v/s\le\sqrt n$. Therefore, $\vartheta\le A$, $\gamma v_i+(Bv)_i\le C$, and
\begin{equation}\label{eq:chain-bounds}
0\le\gamma\le\frac{4\log(e\sqrt n)}{\ell_1\norm v},
\end{equation}
where
\begin{equation}\label{eq:shift}
\begin{gathered}
A=\frac2{\ell_1}\log(1+n^3e^{H-1}),
\\
C=\frac{(5\sqrt m+4)\log(e\sqrt n)}{\ell_1}.
\end{gathered}
\end{equation}
The derivation of these bounds is given in Appendix~\ref{app:lower-chain}. The smallness threshold equals $\delta=1/(4\ell_1\sqrt m)$, and the shift is $\sigma_i=A+2C+4(m-\pos(i))\delta$.

\proofdown{lem:curvature}
\begin{lemma}[membership in the class]\label{lem:curvature}
For each $u\in\R^n$, inclusion~\eqref{eq:inverse-relation} has a unique
solution $v=\Phi(u)\in K$. The operator $\Phi$ is monotone, satisfies
condition~\eqref{eq:smooth} with constants $\ell_0$, $\ell_1$,
$\norm{\Phi(u)}\le G_1$ for all $u$, and $\Phi(-\sigma)=0$.
\end{lemma}

\proofdown{lem:chain}
\begin{lemma}[revealing one coordinate at a time]\label{lem:chain}
For any $u\in\R^n$, the following properties hold.
\begin{enumerate}
\item[(a)] If $|u_n|\le\delta$, then $\sum_i\Phi_i(u)/G_i=1$ and
$\norm{\Phi(u)}\ge G_n/\sqrt n$.
\item[(b)] If $|u_i|\le\delta$ for all $i\ge j$, then $\Phi_i(u)=0$
for $i>j$, and $\Phi(u)$ depends only on $u_1,\ldots,u_j$. In particular,
if $|u_i|\le\delta$ for all $i$, then $\Phi(u)=G_1e_1$.
\end{enumerate}
\end{lemma}

\subsection{How the radius determines the chain length}
\label{sec:lower-radius}
The zero of the constructed operator is located at $u^*=-\sigma$: there $v=0$ solves the inclusion, since $Q\ge Q(0)=0$. Hence, we need to fit $\sigma$ into the ball of radius $R$.

Each coordinate of the shift is of order $\ell_1\sigma_i=\Theta(H+\sqrt m\log n)$.
The term $H$ comes from $\log G_1$: the first answers are exponentially large in the number of blocks. There are $n$ coordinates, so
\[
\ell_1\norm\sigma=\Theta\Bigl(
\underbrace{H\sqrt n}_{\text{many blocks}}+
\underbrace{\frac{n\log n}{\sqrt H}}_{\text{many coordinates per block}}\Bigr).
\]

We choose $m=h^2$, $H=\lceil h\log(4h)\rceil$, and the largest admissible $h$. Then
$n=\Theta(h^3\log h)$ and $\ell_1R=\Theta(h^{5/2}\log^{3/2}h)$. Hence $(\ell_1R)^{6/5}=\Theta(h^3\log^{9/5}h)$, so
\[
n=\Theta\left(\frac{(\ell_1R)^{6/5}}{\log^{4/5}(2+\ell_1R)}\right).
\]
Next, we make the construction smooth and bound the useful radius of the queries.

\subsection{Smoothing preserves the chain}
\label{sec:lower-smooth}
The resulting $\Phi$ is Lipschitz but not necessarily twice differentiable. We take a very small ball of radius $\rho$ and average the values over it:
\[
\Phi_\rho(u)=\int\Phi(u-y)k_\rho(y)\,dy,
\]
where $k_\rho\ge0$, the integral of the kernel equals $1$, and the kernel itself is smooth.

\proofdown{lem:smoothing}
\begin{lemma}[smoothing]\label{lem:smoothing}
Let the support of the kernel $k_\rho\in C^\infty$ lie in the ball of radius
\begin{equation}\label{eq:smoothing-radius}
 0<\rho\le\min\Bigl\{\frac\delta2,\
 \frac{\ell_0}{2\ell_1(\ell_0+\ell_1G_1)},\
 \frac{G_n}{2\ell_1(G_1+G_n)},\ \frac R{16\sqrt n}\Bigr\}.
\end{equation}
Then $\Phi_\rho\in C^\infty$ is monotone, satisfies
condition~\eqref{eq:smooth} with constants $2\ell_0$, $\ell_1$,
$\norm{\Phi_\rho}\le G_1$, and $\Phi_\rho(-\sigma-\rho\one)=0$. The statements
of Lemma~\ref{lem:chain} hold for $\Phi_\rho$ with the threshold $\theta=\delta/2$
instead of $\delta$.
\end{lemma}

The face of the simplex is given by a linear equation, so an average of answers lying on the face also lies on it. Averaging the norm inequality itself would not be possible: large vectors can cancel out.

\subsection{Coordinate hiding}\label{sec:hiding}

\subsubsection{Hiding the coordinates}
\label{sec:lower-local}
If we hide the coordinates without additional protection, the algorithm may query a vector \emph{all} of whose coordinates are greater than $\theta$. Then the chain condition does not apply.

Coordinates were hidden under a bounded norm of the queries by \citet{woodworth2016finite} and \citet{carmon2020lower}. Here we additionally need to preserve $(L_0,L_1)$-Lipschitz continuity with constants that do not depend on the magnitude of the operator. To this end, we multiply the operator by a smooth approximation of the indicator of the segment $[0,1]$, for example,
\[
\chi(t)=
\begin{cases}
1,&t\le1,\\
1-10(t-1)^3+15(t-1)^4-6(t-1)^5,&1<t<2,\\
0,&t\ge2.
\end{cases}
\]
This function is twice continuously differentiable, $0\le\chi\le1$, and $|\chi'|\le15/8$.

\proofdown{lem:localization}
\begin{lemma}\label{lem:localization}
Let $F$ satisfy Assumptions~\ref{asm:monotone}--\ref{asm:root}
with $z_0=0$ and $\norm{F(z)}\le M$ for all $z$. Let $L=\max\{4R,32M/\mu\}$ and
\[
 0<\mu\le\frac{L_0}{8(1+L_1R)}.
\]
Then the operator
\begin{equation}\label{eq:localization}
 \widetilde F(z)=\chi\bigl(\norm z/L\bigr)F(z)+\mu z
\end{equation}
is $\mu/2$-strongly monotone, has a unique zero $z^\ast_\mu$ with
$\norm{z^\ast_\mu}\le\norm{z^\ast}$, and satisfies
condition~\eqref{eq:smooth} with constants $2L_0$, $2L_1$.
For $\norm z\le L$ it equals $F(z)+\mu z$, and for $\norm z\ge2L$ it equals the known
operator $\mu z$. If $F\in C^2$, then also $\widetilde F\in C^2$.
\end{lemma}

\proofdown{lem:hardness}
\begin{lemma}[preservation of the success condition]\label{lem:hardness}
Let the conditions of Lemma~\ref{lem:localization} hold.
\begin{enumerate}
\item[(a)] If $\mu\le\varepsilon/(16R)$, $M\ge\varepsilon$, and
$\norm{\widetilde F(z)}\le\varepsilon$, then $\norm z<L$ and $\norm{F(z)}\le2\varepsilon$.
\item[(b)] If $\norm{z^\ast}\le R/4$, $\mu\le\eta/(16R^2)$, $M\ge\eta/R$, and
$\Gap^{\widetilde F}_R(z)\le\eta$, then $\norm z\le L/4$ and $\Gap^F_R(z)<4\eta$.
\end{enumerate}
\end{lemma}

\subsubsection{Random embedding}
\label{sec:lower-embed}
We choose in advance random distinct positions $\pi(1),\ldots,\pi(n)$ among the $d$ coordinates. The matrix $U=(e_{\pi(1)},\ldots,e_{\pi(n)})$ simply places the coordinates of the short chain at these positions: $(U^\top z)_i=z_{\pi(i)}$. The full operator has the form $F_U(z)=\chi(\norm z/L)U\Phi_\rho(U^\top z)+\mu z$.
The embedding preserves the norms of the chain vectors and monotonicity. For $\norm z\ge2L$, the oracle answer $\mu z$ does not depend on $U$ at all.

\proofdown{prop:random-hiding}
\begin{proposition}[random embedding]\label{prop:random-hiding}
Let $n\ge2$, let the operator $\Phi\colon\R^n\to\R^n$ satisfy
Assumptions~\ref{asm:monotone}--\ref{asm:root} with $z_0=0$ and
$\norm{\Phi(u)}\le M$ for all $u$, and let the numbers $\mu$ and $L$ satisfy
the conditions of Lemmas~\ref{lem:localization} and~\ref{lem:hardness}.
Suppose that for some $\theta>0$
\begin{enumerate}
\item[(a)] if $|u_i|\le\theta$ for all $i\ge j$, then $\Phi_i(u)=0$ for $i>j$,
and $\Phi(u)$ depends only on $u_1,\ldots,u_j$,
\item[(b)] if $|u_n|\le\theta$, then $\norm{\Phi(u)}>2\varepsilon$
(respectively, $\Gap^\Phi_R(u)>4\eta$).
\end{enumerate}
Let $N=\lfloor n/2\rfloor$, $d\ge n+48(N+1)nL^2/\theta^2$,
and let $U\in\R^{d\times n}$ be a random matrix whose columns are distinct
standard basis vectors of $\R^d$ chosen uniformly.
Then the operators $F_U(z)=\chi\bigl(\norm z/L\bigr)\,U\Phi(U^\top z)+\mu z$ belong to $\mathcal P_d$ with constants $2L_0$, $2L_1$. If an algorithm,
on every operator $F_U$, outputs with probability at least $2/3$ a point $\hat z$
with $\norm{F_U(\hat z)}\le\varepsilon$ (respectively, with
$\Gap_R^{F_U}(\hat z)\le\eta$), then on one of the operators $F_U$
it makes at least $7(N+1)/12$ oracle calls in expectation.
\end{proposition}

\subsection{Choice of parameters and the lower bound for the norm}
\label{sec:lower-assemble}
It remains to verify that, for the required $n$, the error is indeed greater than $\varepsilon$, the zero is close enough, and the $(L_0,L_1)$-Lipschitz constants coincide with the original ones.

For $L_1=0$ there is no nonlinear term, so from now on $L_1>0$. Let $x=L_1R$ and $q=L_0R/\varepsilon\ge1$. We choose $\ell_0=L_0/4$, $\ell_1=L_1/2$, $m=h^2$, $H=\lceil h\log(4h)\rceil$, and $n=mH$, where $h\ge256$ is the largest integer such that
\begin{equation}\label{eq:h-choice}
2048h^{5/2}\log^{3/2}(4h)\le x.
\end{equation}
If no such $h$ exists, then $x$ is bounded by an absolute constant and the nonlinear term is absorbed into the constant term $1$ in $1+\mathcal T^\ast$.

Lemma~\ref{lem:radius} and the maximality of $h$ give $\norm\sigma\le R/8$, $n^{5/2}>4x$, and $n^3>4x^2$. After smoothing, the zero moves by at most $R/16$. Hence, it remains within $R/4$.

When the last coordinate is small, the smoothed chain has norm
\[
\norm{\Phi_\rho(u)}\ge\frac{G_n}{\sqrt n}
=\frac{L_0}{2L_1}n^{5/2}>2L_0R\ge2\varepsilon.
\]
This is the required large error. We choose $L=\max\{4R,32G_1/\mu\}$,
\[
\mu=\frac{L_0}{64(1+q+x)},
\]
and the dimension from Proposition~\ref{prop:random-hiding}. The conditions of Lemmas~\ref{lem:localization} and~\ref{lem:hardness} hold, and the verification of the inequalities is given in Appendix~\ref{app:lower-radius}.

The $(L_0,L_1)$-Lipschitz constants go through three stages:
\[
(\ell_0,\ell_1)
\xrightarrow{\text{smoothing}}(2\ell_0,\ell_1)
\xrightarrow{\text{coordinate hiding}}(4\ell_0,2\ell_1)=(L_0,L_1).
\]
We have obtained admissible $C^2$ operators. For every algorithm, among them there is a worst-case operator with $T\gtrsim n$. Therefore,
\[
1+\mathcal T^\ast(\varepsilon)\gtrsim
1+\frac{(L_1R)^{6/5}}{\log^{4/5}(2+L_1R)}.
\]
The exponent $6/5$ arises entirely from the radius constraint. Smoothing, hiding, and the probabilistic reduction change only absolute factors and the dimension.

\subsection{Passing to the Minty gap}
\label{sec:lower-gap}
The gap equals $\Gap_R(u)=\sup_{\norm w\le R}\ip{\Phi(w)}{u-w}$. To bound it from below, it suffices to substitute \emph{one} suitable point $w$ into the supremum.

We choose a query $w$ whose answer is the last vertex of the simplex, $G_ne_n$. Then the inner product becomes one-dimensional: $\ip{\Phi(w)}{u-w}=G_n(u_n-w_n)$. It only remains to place $w_n$ noticeably to the left of zero while $u_n$ stays small. The definition of the inverse operator directly suggests the required query: to $-\sigma+\nabla Q(G_ne_n)+G_nBe_n$ we add the outer normal $e_n/\ell_1$ at the vertex $G_ne_n$.

\proofdown{lem:witness}
\begin{lemma}[test point for the gap]\label{lem:witness}
Let
\[
 w=-\sigma+\nabla Q(G_ne_n)+G_nBe_n+\frac{e_n}{\ell_1},
\]
and let $\rho$ satisfy~\eqref{eq:smoothing-radius}. Then
$\Phi_\rho(w)=G_ne_n$, $w_n\le-3/\ell_1$, and
$\ip{\Phi_\rho(w)}{u-w}\ge G_n/\ell_1$ for $|u_n|\le\delta$.
With the parameters of Section~\ref{sec:lower-assemble}, also $\norm w\le R/8$.
\end{lemma}

Hence, for small $u_n$,
\[
\Gap_R^{\Phi_\rho}(u)\ge\frac{G_n}{\ell_1}
=\frac{L_0}{L_1^2}n^3>4L_0R^2\ge4\eta.
\]
Here we require $L_0R^2/\eta\ge1$.

In the choice of $\mu$, we replace $q$ by $L_0R^2/\eta$. Then the conditions of Lemma~\ref{lem:hardness}(b) hold, and the same random embedding gives the same cost of order $n$ for the gap.

\subsection{The final bound}
\label{sec:linear-lower}
Lipschitz operators belong to the class under consideration for any $L_1\ge0$. For deterministic methods, lower bounds have long been known~\citep{nemirovsky1992}. For our randomized model, we need the following fact: for symmetric $A$ with spectrum in
$[-\kappa,-1]\cup[1,\kappa]$, an arbitrary randomized solver with relative residual $1/2$ requires $1+T\gtrsim\kappa$ matrix--vector multiplications in expectation.

\proofdown{lem:linear-budget}
\begin{lemma}[solving linear systems]\label{lem:linear-budget}
Suppose that an algorithm using arbitrary adaptive multiplications $v\mapsto Av$,
for every symmetric matrix $A$ with spectrum in $[-\kappa,-1]\cup[1,\kappa]$
and every vector $b\ne0$, almost surely stops and outputs $\hat x$ with
\[
 \Pr\bigl\{\norm{A\hat x-b}\le\norm b/2\bigr\}\ge\frac23.
\]
If $T$ is the largest expected number of multiplications over $A,b$ and the dimension,
then $1+T\gtrsim\kappa$.
\end{lemma}

The lemma follows from the result of \citet[Theorems~2.6--2.8]{derezinski2026}. The passage to success probability $2/3$ and to a bound on the expected number of multiplications is proved in Appendix~\ref{app:lower-linear}.

We embed the system into a monotone problem. Take
\[
F(x,y)=a\left(Ay,\frac{s b}{\norm b}-Ax\right).
\]
The two blocks with opposite signs give a skew-symmetric linear part. Therefore, $F$ is monotone, $\norm{DF}=a\norm A\le a\kappa$, and its zero
\[
z^*=\left(\frac{sA^{-1}b}{\norm b},0\right)
\]
satisfies $\norm{z^*}\le s$. Each query of $F$ costs two multiplications by $A$.

For the norm, we choose $a=2\varepsilon/R$, $s=R$, and $\kappa=L_0R/(2\varepsilon)$. If $\norm{F(x,y)}\le\varepsilon$, then
\[
\left\|A\frac{\norm b}{R}x-b\right\|
\le\frac{\norm b}{R}\frac{\varepsilon}{a}=\frac{\norm b}{2}.
\]
Consequently, $1+\mathcal T^\ast(\varepsilon)\gtrsim L_0R/\varepsilon$.

For the gap, we take $a=16\eta/R^2$, $s=R/4$, and $\kappa=L_0R^2/(16\eta)$. For an affine skew-symmetric operator,
\[
\Gap_R^F(z)=R\norm{F(z)}+\ip{F(z)}{z^*}
\ge\frac{3R}{4}\norm{F(z)}.
\]
The condition $\Gap_R^F\le\eta$ gives a relative residual of the system of at most $1/3$, so $1+\mathcal T^\ast_{\Gap}(\eta)\gtrsim L_0R^2/\eta$. The details of the identity for the gap are given in Appendix~\ref{app:lower-linear}.

Two different problems give two lower bounds, and their maximum is at least half of their sum. Therefore, the linear problem and the chain together give the claimed bound. When these ratios are bounded, the corresponding lower bounds follow from $1+\mathcal T^\ast\ge1$.

\proofdown{prop:linear-lower}
\begin{proposition}[the term $L_0R/\varepsilon$]\label{prop:linear-lower}
If $L_0R/\varepsilon\ge1$, then
\[
 1+\mathcal T^\ast(\varepsilon)\gtrsim\frac{L_0R}{\varepsilon}.
\]
If $L_0R^2/\eta\ge1$, then
\[
 1+\mathcal T^\ast_{\Gap}(\eta)\gtrsim\frac{L_0R^2}{\eta}.
\]
\end{proposition}

\begin{theorem}[lower bound with exponent $6/5$]\label{thm:nonlinear-lower}
For $L_0R/\varepsilon\ge1$,
\[
 1+\mathcal T^\ast(\varepsilon)\ge c\Bigl(1+\frac{L_0R}{\varepsilon}
 +\frac{(L_1R)^{6/5}}{\log^{4/5}(2+L_1R)}\Bigr),
\]
and the same bound with $L_0R/\varepsilon$ replaced by $L_0R^2/\eta$ holds for
$\mathcal T^\ast_{\Gap}(\eta)$ when $L_0R^2/\eta\ge1$. The worst-case operators can
be chosen from $C^2(\R^d)$.
\end{theorem}

\begin{proof}[Completion of the proof of Theorem~\ref{thm:nonlinear-lower}]
In Sections~\ref{sec:lower-assemble} and~\ref{sec:lower-gap}, we constructed admissible $C^2$ operators for which Proposition~\ref{prop:random-hiding} gives a lower bound of order $n$. By Lemma~\ref{lem:radius} and the maximality of $h$, this is of order $(L_1R)^{6/5}/\log^{4/5}(2+L_1R)$. Proposition~\ref{prop:linear-lower} independently gives the Lipschitz term. Hence, there exist absolute constants $c_1,c_2>0$ for which
\[
1+T\ge c_1\frac{L_0R}{\varepsilon},
\]
\[
1+T\ge c_2\frac{(L_1R)^{6/5}}{\log^{4/5}(2+L_1R)}.
\]
Adding these inequalities to $1+T\ge1$, we obtain
\[
1+T\ge\frac{\min\{1,c_1,c_2\}}3\Bigl(1+\frac{L_0R}{\varepsilon}+\frac{(L_1R)^{6/5}}{\log^{4/5}(2+L_1R)}\Bigr).
\]
For the gap, the argument is the same with $L_0R/\varepsilon$ replaced by $L_0R^2/\eta$. When $L_1R$ is bounded, including $L_1=0$, the nonlinear term is bounded by an absolute constant and is absorbed into the constant term $1$. Taking the infimum over admissible algorithms, we obtain both statements.
\end{proof}

The exponent $6/5$ and the logarithmic factor arise from the relation between the chain length $n$ and the distance $\norm\sigma$ from the starting point to the solution: with $n$ coordinates, the shift $\sigma$ that defines the solution $u^*=-\sigma$ has norm of order $\sqrt n(H+\sqrt m\log n)/\ell_1$, and this norm must not exceed $R$.
Smoothing and coordinate hiding change only absolute constants and the dimension, without reducing the required number of queries. The initial norm of the worst-case operator equals $G_1=(L_0/(2L_1))n^3e^{H-1}$ and can be exponentially large in the number of blocks.

\section{Variational inequalities on a bounded set}\label{sec:constraints}

In this section, $Z\subset\R^d$ is a convex compact set of diameter $D$ onto which
the exact projection $P_Z$ is available, and the output point must belong to $Z$.
A solution of the variational inequality is a point $z^\ast\in Z$ such that
$\ip{F(z^\ast)}{z-z^\ast}\ge0$ for all $z\in Z$. In other words,
$-F(z^\ast)$ lies in the normal cone
\[
 N_Z(z^\ast)=\{\xi\in\R^d:\ip\xi{z-z^\ast}\le0\ \text{for all}\ z\in Z\},
\]
so $F(z^\ast)$ need not be zero. Accuracy is measured by the Minty gap
and the strong gap
\begin{gather*}
 \Gap^{\rm M}_Z(z)=\sup_{u\in Z}\ip{F(u)}{z-u},\\
 \Gap^{\rm S}_Z(z)=\sup_{u\in Z}\ip{F(z)}{z-u}.
\end{gather*}
For a monotone operator, $0\le\Gap^{\rm M}_Z\le\Gap^{\rm S}_Z$, and the equality
$\Gap^{\rm M}_Z(z)=0$ means that $z$ is a solution \citep{minty1962}.
The relation between these gaps and the residuals used by projection methods
is given in Appendix~\ref{app:constraints}.

\subsection{Why the parameters \texorpdfstring{$L_0$, $L_1$, and $D$}{L0, L1, and D} are not sufficient}
On the whole space, a solution requires $F(z^\ast)=0$. On a compact set, the feasible directions may lie only on one side of a boundary point, so a nonzero component of $F(z^\ast)$ opposite to the outward normal may persist at the solution.
Its inner product with any vector $z-z^\ast$, $z\in Z$, is nonnegative.
 For each operator this component is finite, but its magnitude is not bounded in terms of $L_0$, $L_1$, and $D$. If it is large, condition~\eqref{eq:smooth} admits a Jacobian of large norm with the same $L_0,L_1$.

Fix $L_0,L_1,D,\eta$ with $L_1>0$. We introduce a new parameter of the construction, an arbitrary real number $\Gamma>0$, which we will regard as arbitrarily large. We add the constant component $\Gamma/L_1$ to the worst-case linear problem of Section~\ref{sec:linear-lower}. Then $\norm F\ge\Gamma/L_1$ implies $L_0+L_1\norm F\ge\Gamma$.
Therefore, condition~\eqref{eq:smooth} admits a Jacobian of norm $\Gamma$ with the same $L_0,L_1$.

Choose $\Gamma\ge72\eta/D^2$, set $\kappa=\Gamma D^2/(72\eta)$, and take a symmetric matrix $A$ with spectrum in $[-\kappa,-1]\cup[1,\kappa]$ and $b\ne0$ from Lemma~\ref{lem:linear-budget}. On the set $Z=B(0,D/3)\times[0,D/3]\times B(0,D/3)$ of diameter $D$, we construct
\[
F(x,t,y)=\left(\frac{72\eta}{D^2}Ay,\frac\Gamma{L_1},
\frac{72\eta}{D^2}\left(\frac{Db}{12\norm b}-Ax\right)\right).
\]

The exact verification in Appendix~\ref{app:constraints-lower} reduces a small gap to solving a linear system and yields the lower bound $T\gtrsim\Gamma D^2/\eta$. Since $\Gamma$ can be chosen arbitrarily large, we obtain the following result.

\proofdown{thm:constrained-impossible}
\begin{theorem}[no uniform bound]\label{thm:constrained-impossible}
Let $L_0\ge0$ and $L_1,D,\eta>0$ be fixed. For any $\Gamma\ge72\eta/D^2$, there exist
a dimension $d$, a set $Z$ of diameter $D$ that is the product of two balls and a segment, and
a family of affine monotone operators $F$ satisfying
condition~\eqref{eq:smooth} on the whole $\R^d$ such that any randomized
algorithm that outputs, with probability at least $2/3$, a point $\hat z\in Z$
with $\Gap^{\rm M}_Z(\hat z)\le\eta$ makes at least
$c\,\Gamma D^2/\eta$ oracle calls in expectation. The same holds for
$\Gap^{\rm S}_Z$. In particular, the minimax complexity of the class
specified only by the numbers $L_0$, $L_1$, and $D$ is infinite.
\end{theorem}

Thus, on a bounded set, a uniform complexity bound
requires an additional parameter, for example an upper bound on $\norm{F(z^\ast)}$.

\section{Conclusion}\label{sec:conclusion}

For finding a zero of a monotone $(L_0,L_1)$-Lipschitz operator with a known
radius $R$, we determined the oracle complexity up to the factor
$\log^{4/5}(2+L_1R)$: it is of order
$1+L_0R/\varepsilon+(L_1R)^{6/5}$. The upper bound is deterministic, and the lower
bound holds for arbitrary randomized methods. For gradients of
convex functions, the price of traversing the region of large operator values is
$L_1R$ \citep{rudomanov2026everest}, while for monotone operators it
grows to $(L_1R)^{6/5}$.

\subsection{Open questions}\label{sec:open-questions}
\begin{enumerate}
\item \emph{Saddle-point problems.}
What is the complexity for operators of saddle-point problems
$F(x,y)=\bigl(\nabla_x\varphi(x,y),-\nabla_y\varphi(x,y)\bigr)$
with a convex--concave function $\varphi$? For this class, \method{ABSEIL} yields a complexity bound with exponent
$6/5$ in $L_1R$. The operator in our lower bound is not constructed as the operator
of a saddle-point problem, while minimization of a convex function is a special case of a saddle-point
problem for which the price of traversing the region of large values equals
$L_1R$ \citep{rudomanov2026everest}. Thus, for saddle-point problems
the optimal exponent lies between $1$ and $6/5$, and its exact value
remains open.

\item \emph{Unknown parameters.}
Can the same order be achieved by a method that does not know
$L_0$, $L_1$, and $R$, counting all oracle calls spent on stepsize tuning
and on checking the stopping criterion? The steps of our method are adaptive, but they use
the known $L_0$, $L_1$, and $R$.
\end{enumerate}

\begingroup
\renewcommand{\refname}{References}
\renewcommand{\bibfont}{\small}
\bibliographystyle{plainnat}
\bibliography{references}
\endgroup
\clearpage
\appendix
\begingroup
\small
\section{Properties of the class}\label{app:tools}

\subsection{Characterization of the class}

\begin{proof}[Proof of Lemma~\ref{lem:class}\proofup{lem:class}]
\leavevmode
\begin{itemize}
\item \emph{Necessity in~(a).} A continuous $F$ is bounded on each ball
by some number $M$. If the segment $[x,y]$ lies in such a ball, then
by~\eqref{eq:smooth} $\norm{F(x)-F(y)}\le(L_0+L_1M)\norm{x-y}$.
Hence, $t\mapsto F(x+tv)$ is locally Lipschitz, locally absolutely continuous,
and differentiable for almost all $t$. At such a point $t$, we apply~\eqref{eq:smooth}
to the points $x+tv$ and $x+(t+h)v$, divide by $|h|$, and let $h$ tend to zero.
By continuity of $F$, the supremum in~\eqref{eq:smooth} tends to $\norm{F(x+tv)}$.
\item \emph{Sufficiency in~(a).} Let $v=y-x$. By the Newton--Leibniz formula,
\begin{align*}
 \norm{F(y)-F(x)}
 &\le\int_0^1\Bigl\|\frac{d}{dt}F(x+tv)\Bigr\|\,dt
 \le\int_0^1\bigl(L_0+L_1\norm{F(x+tv)}\bigr)\norm v\,dt\\
 &\le\Bigl(L_0+L_1\sup_{w\in[x,y]}\norm{F(w)}\Bigr)\norm{x-y}.
\end{align*}
\item \emph{Claim~(b).} For $F\in C^1$, the derivative of the map
$t\mapsto F(x+tv)$ equals $DF(x+tv)v$ and is continuous in $t$. Therefore, condition~(a)
is equivalent to the inequality $\norm{DF(z)v}\le(L_0+L_1\norm{F(z)})\norm v$ for all
$z,v$, that is, to the bound $\norm{DF(z)}\le L_0+L_1\norm{F(z)}$.
\item \emph{Claim~(c).} By Rademacher's theorem, $F$ is differentiable almost
everywhere. Fix a direction $v$. By Fubini's theorem, for almost all $x$
in the orthogonal complement of $v$, the derivative $DF(x+tv)$ exists
for almost all $t$ and equals $\frac d{dt}F(x+tv)$. On such lines,
condition~(a) holds, and by what has been proved,~\eqref{eq:smooth} holds for segments
parallel to $v$ and lying on these lines. An arbitrary segment with direction
$v$ is approximated by such segments, and~\eqref{eq:smooth} passes to the limit
by continuity of $F$.\qedhere
\end{itemize}
\end{proof}

\subsection{Grönwall's lemma}

\begin{proof}[Proof of Lemma~\ref{lem:gronwall}\proofup{lem:gronwall}]
Let $s=\norm{u-v}>0$, $w(t)=u+t(v-u)/s$, and $\phi(t)=L_0+L_1\norm{F(w(t))}$
for $t\in[0,s]$. By Lemma~\ref{lem:class}(a), the function $\phi$ is absolutely
continuous, and for almost all $t$,
\[
 \phi'(t)\le L_1\Bigl\|\frac d{dt}F(w(t))\Bigr\|\le L_1\phi(t).
\]
Therefore, $\bigl(e^{-L_1t}\phi(t)\bigr)'\le0$ and $\phi(t)\le\phi(0)e^{L_1t}$,
which for $t=s$ gives~(a). Further,
\[
 \norm{F(v)-F(u)}
 \le\int_0^s\Bigl\|\frac d{dt}F(w(t))\Bigr\|dt
 \le\int_0^s\phi(t)\,dt
 \le\phi(0)\int_0^se^{L_1t}dt
 =\phi(0)\frac{e^{L_1s}-1}{L_1}
 \by{\le}{(\mathrm a)}\phi(0)\,s\,e^{L_1s},
\]
where (a) is the inequality $e^x-1\le xe^x$. If the derivative bound holds
only on the segment $[u,v]$ and with a constant $A$ in place of $L_0$, then the same argument
goes through with $\phi(t)=A+L_1\norm{F(w(t))}$.
\end{proof}

\subsection{Gap and the norm of \texorpdfstring{$F$}{F}}

\begin{proof}[Proof of Lemma~\ref{lem:gap-norm}\proofup{lem:gap-norm}]
For $u=z^\ast$, we have $\ip{F(z^\ast)}{z-z^\ast}=0$, hence $\Gap_R(z)\ge0$.
For any $u\in B(z_0,R)$,
\[
 \ip{F(u)}{z-u}\by{\le}{(\mathrm a)}\ip{F(z)}{z-u}
 \le\norm{F(z)}\norm{z-u}
 \by{\le}{(\mathrm b)}\norm{F(z)}\bigl(\norm{z-z_0}+R\bigr),
\]
where (a) is monotonicity and (b) is the inequality $\norm{u-z_0}\le R$.
\end{proof}

\section{Proofs of the upper bound}\label{app:upper}

Throughout this appendix, Assumptions~\ref{asm:monotone}--\ref{asm:root} hold.

\subsection{Long step}

\begin{proof}[Proof of Lemma~\ref{lem:long-step}\proofup{lem:long-step}]
We prove the lemma for an arbitrary $b\in[1,L_1R]$, in the form used
in Section~\ref{sec:total}. Assume that $\norm{F(x)}>2L_0R(e^b-1)/b$.
We need to show that a step of length $b/(2L_1)$ decreases the squared distance
to the solution by $b^2/(4L_1^2)$. Lemma~\ref{lem:long-step} is obtained for
$b=(L_1R)^{2/5}$.

Consider the point $w=z^\ast+bF(x)/(L_1\norm{F(x)})$. Since
$\norm{w-z^\ast}=b/L_1$, Lemma~\ref{lem:gronwall}(b) with $u=z^\ast$
gives $\norm{F(w)}\le(L_0/L_1)(e^b-1)$. We show that
$\ip{F(x)}{x-z^\ast}>b\norm{F(x)}/(2L_1)$. Suppose the contrary. Then
\[
 \frac b{L_1}-\frac{\ip{F(x)}{x-z^\ast}}{\norm{F(x)}}\ge\frac b{2L_1}>0
\]
and
\begin{align*}
 \norm{x-w}^2
 &=\norm{x-z^\ast}^2-\frac{2b}{L_1}\frac{\ip{F(x)}{x-z^\ast}}{\norm{F(x)}}+\frac{b^2}{L_1^2}\\
 &\by{\le}{(\mathrm a)}R^2-\frac{b^2}{L_1^2}
 +\frac{2b}{L_1}\Bigl(\frac b{L_1}-\frac{\ip{F(x)}{x-z^\ast}}{\norm{F(x)}}\Bigr)\\
 &\by{\le}{(\mathrm b)}\frac{4R^2L_1^2}{b^2}
 \Bigl(\frac b{L_1}-\frac{\ip{F(x)}{x-z^\ast}}{\norm{F(x)}}\Bigr)^2,
\end{align*}
where (a) is the condition $\norm{x-z^\ast}\le R$, and (b) is the inequality
$R^2-\rho^2+2\rho\sigma\le4R^2\sigma^2/\rho^2$, which holds for $\sigma\ge\rho/2$
and $\rho\le R$, since the difference of its sides equals
$(2\sigma-\rho)(2R^2\sigma+R^2\rho-\rho^3)/\rho^2\ge0$. It is applied
with $\rho=b/L_1$ and $\sigma=b/L_1-\ip{F(x)}{x-z^\ast}/\norm{F(x)}$.
By monotonicity and the Cauchy--Schwarz inequality,
\begin{align*}
 \norm{F(x)}\Bigl(\frac b{L_1}-\frac{\ip{F(x)}{x-z^\ast}}{\norm{F(x)}}\Bigr)
 &=\ip{F(x)}{w-x}\le\ip{F(w)}{w-x}\le\norm{F(w)}\norm{x-w}\\
 &\le\frac{L_0}{L_1}(e^b-1)\frac{2RL_1}{b}
 \Bigl(\frac b{L_1}-\frac{\ip{F(x)}{x-z^\ast}}{\norm{F(x)}}\Bigr).
\end{align*}
Dividing by the positive factor, we obtain
$\norm{F(x)}\le2L_0R(e^b-1)/b$, which contradicts the assumption. Consequently,
\[
 \Bigl\|x-\frac b{2L_1}\frac{F(x)}{\norm{F(x)}}-z^\ast\Bigr\|^2
 =\norm{x-z^\ast}^2-\frac b{L_1}\frac{\ip{F(x)}{x-z^\ast}}{\norm{F(x)}}+\frac{b^2}{4L_1^2}
 <\norm{x-z^\ast}^2-\frac{b^2}{4L_1^2}.\qedhere
\]
\end{proof}

\subsection{Two cumulative bounds}

\begin{proof}[Proof of Lemma~\ref{lem:energies}\proofup{lem:energies}]
Denote by
\[
 \mathcal L_k=\ip{u_k-y_k}{F(y_k)}-\frac{\Lambda_k}3\norm{F(y_k)}^2
\]
the Lyapunov function from \citet[Proposition~2]{sun2026}.
For $k=1$, the factor $1-\alpha_1$ equals zero, and the value $F(y_0)$ does not enter
any formula. From the definition of Algorithm~\ref{alg:belay}, we obtain
\begin{equation}\label{eq:maeg-identities}
\begin{aligned}
 &y_k-y_{k-1}=\alpha_k(u_{k-1}-y_{k-1})-\lambda_kF(x_k),\\
 &u_k-y_k=(1-\alpha_k)(u_{k-1}-y_{k-1})+\lambda_k\bigl(F(x_k)-F(y_k)\bigr)+\tfrac56\lambda_kF(y_k),\\
 &x_k-y_k=\lambda_k\bigl(F(x_k)-(1-\alpha_k)F(y_{k-1})\bigr),\\
 &\frac{\lambda_k(1-\alpha_k)}{\alpha_k}=\frac{2\Lambda_{k-1}}3.
\end{aligned}
\end{equation}
\begin{itemize}
\item \emph{Lyapunov function.} Substituting~\eqref{eq:maeg-identities}
into the definition of $\mathcal L_k$, by direct computation we obtain the identity
\begin{align*}
 \mathcal L_k-(1-\alpha_k)\mathcal L_{k-1}
 &=\frac{\lambda_k}{2\alpha_k}\Bigl(\norm{F(x_k)-(1-\alpha_k)F(y_{k-1})}^2-\norm{F(x_k)-F(y_k)}^2\Bigr)\\
 &\quad+\frac{1-\alpha_k}{\alpha_k}\ip{F(y_k)-F(y_{k-1})}{y_k-y_{k-1}}.
\end{align*}
The second term is nonnegative by monotonicity, and by~\eqref{eq:accept}
and~\eqref{eq:maeg-identities},
\[
 \norm{F(x_k)-(1-\alpha_k)F(y_{k-1})}=\frac{\norm{x_k-y_k}}{\lambda_k}
 \ge4\norm{F(x_k)-F(y_k)}.
\]
Consequently, $\mathcal L_k\ge(1-\alpha_k)\mathcal L_{k-1}$, and $\mathcal L_1\ge0$
since $\alpha_1=1$. By induction, $\mathcal L_k\ge0$ for all $k$.
\item \emph{Distance from the anchor to the solution.} We have
\[
 \ip{u_k-z^\ast}{F(y_k)}
 =\mathcal L_k+\frac{\Lambda_k}3\norm{F(y_k)}^2+\ip{F(y_k)-F(z^\ast)}{y_k-z^\ast}
 \by{\ge}{(\mathrm a)}\frac{\Lambda_k}3\norm{F(y_k)}^2,
\]
where (a) is the inequality $\mathcal L_k\ge0$ together with monotonicity. Since
$u_{k-1}=u_k+\lambda_kF(y_k)/6$,
\begin{align*}
 \norm{u_k-z^\ast}^2
 &=\norm{u_{k-1}-z^\ast}^2-\frac{\lambda_k}3\ip{u_k-z^\ast}{F(y_k)}
 -\frac{\lambda_k^2}{36}\norm{F(y_k)}^2\\
 &\le\norm{u_{k-1}-z^\ast}^2-\frac{\lambda_k\Lambda_k}9\norm{F(y_k)}^2.
\end{align*}
Summing and using $\norm{u_0-z^\ast}=\norm{s-z^\ast}\le R$, we obtain
\begin{equation}\label{eq:anchor-sum}
 \norm{u_k-z^\ast}^2+\frac19\sum_{i=1}^k\lambda_i\Lambda_i\norm{F(y_i)}^2\le R^2.
\end{equation}
These are the first and third bounds of the lemma. The second one follows from
\[
 R\norm{F(y_k)}\ge\norm{u_k-z^\ast}\norm{F(y_k)}\ge\ip{u_k-z^\ast}{F(y_k)}
 \ge\frac{\Lambda_k}3\norm{F(y_k)}^2.
\]
\item \emph{Errors of the extragradient step.} By~\eqref{eq:accept}
and~\eqref{eq:maeg-identities},
\begin{align*}
 4\lambda_k\norm{F(x_k)-F(y_k)}
 &\le\norm{x_k-y_k}\\
 &\le\lambda_k\bigl(\norm{F(y_k)}+\norm{F(x_k)-F(y_k)}+(1-\alpha_k)\norm{F(y_{k-1})}\bigr),
\end{align*}
whence
\[
 \norm{F(x_k)-F(y_k)}^2\le\frac29\bigl(\norm{F(y_k)}^2+(1-\alpha_k)^2\norm{F(y_{k-1})}^2\bigr).
\]
Keeping the factor $(1-\alpha_k)^2$, for $k\ge2$ we obtain
$\lambda_k^2/\alpha_k\le\lambda_k\Lambda_k$ and
\[
 \frac{\lambda_k^2(1-\alpha_k)^2}{\alpha_k}
 =\frac{4\lambda_k\Lambda_{k-1}^2}{3(2\Lambda_{k-1}+3\lambda_k)}
 \le\frac23\lambda_k\Lambda_{k-1}
 \by{\le}{(\mathrm b)}\frac83\lambda_{k-1}\Lambda_{k-1},
\]
where (b) is the condition $\lambda_k\le4\lambda_{k-1}$. Therefore,
\begin{equation}\label{eq:error-sum}
 \sum_{i=1}^k\frac{\lambda_i^2}{\alpha_i}\norm{F(x_i)-F(y_i)}^2
 \by{\le}{\eqref{eq:anchor-sum}}\frac29\Bigl(9+\frac83\cdot9\Bigr)R^2=\frac{22}3R^2.
\end{equation}
\item \emph{Distance from the point to the anchor.} A direct computation
using~\eqref{eq:maeg-identities} gives
\begin{align*}
 &\bigl(\norm{u_k-y_k}^2+5\norm{u_k-z^\ast}^2\bigr)
 -\bigl(\norm{u_{k-1}-y_{k-1}}^2+5\norm{u_{k-1}-z^\ast}^2\bigr)\\
 &\qquad=\bigl\|(1-\alpha_k)(u_{k-1}-y_{k-1})+\lambda_k\bigl(F(x_k)-F(y_k)\bigr)\bigr\|^2
 -\norm{u_{k-1}-y_{k-1}}^2\\
 &\qquad\quad-\frac53\lambda_k\ip{F(y_k)}{y_k-z^\ast}-\frac56\lambda_k^2\norm{F(y_k)}^2\\
 &\qquad\by{\le}{(\mathrm c)}\frac{\lambda_k^2}{\alpha_k}\norm{F(x_k)-F(y_k)}^2
 -\alpha_k\norm{u_{k-1}-y_{k-1}}^2,
\end{align*}
where (c) is the convexity of the squared norm,
$\norm{(1-\alpha)a+\alpha b}^2\le(1-\alpha)\norm a^2+\alpha\norm b^2$
with $b=\lambda_k(F(x_k)-F(y_k))/\alpha_k$, together with the inequality
$\ip{F(y_k)}{y_k-z^\ast}\ge0$. Summing and using $u_0=y_0=s$
and~\eqref{eq:error-sum}, we obtain
\[
 \norm{u_k-y_k}^2+5\norm{u_k-z^\ast}^2+\sum_{i=1}^k\alpha_i\norm{u_{i-1}-y_{i-1}}^2
 \le5R^2+\frac{22}3R^2=\frac{37}3R^2.\qedhere
\]
\end{itemize}
\end{proof}

\subsection{Number of iterations of the second stage}

\begin{proof}[Proof of Proposition~\ref{prop:adaptive}\proofup{prop:adaptive}]
By construction, $\lambda_1(L_0+L_1\norm{F(s)})=1/6$ and $\lambda_k(L_0+L_1B_k)\le1/6$
for $k\ge2$. We use the inequalities $e^{1/6}<6/5$ and $e^{1/5}<5/4$.
\begin{itemize}
\item \emph{Admissibility of the step.} For $k=1$, we have $x_1=s$ and
$\norm{x_1-y_1}=\lambda_1\norm{F(s)}$, hence by Lemma~\ref{lem:gronwall},
\[
 \lambda_1\norm{F(x_1)-F(y_1)}\le\frac{e^{1/6}}6\norm{x_1-y_1}<\frac15\norm{x_1-y_1}.
\]
For $k\ge2$, we have $\alpha_k\le3\lambda_k/(2\Lambda_{k-1})$, hence the equality
$x_k-y_{k-1}=\alpha_k(u_{k-1}-y_{k-1})-(1-\alpha_k)\lambda_kF(y_{k-1})$
implies $\norm{x_k-y_{k-1}}\le\lambda_kB_k$, and $L_1\lambda_kB_k\le1/6$. Applying
Lemma~\ref{lem:gronwall} twice, we obtain
\begin{gather*}
 \norm{F(x_k)-F(y_{k-1})}\le\frac{e^{1/6}}6B_k<\frac15B_k,\\
 \norm{F(x_k)}\le\frac65B_k,\\
 \norm{x_k-y_k}\le\frac65\lambda_kB_k,\\
 \lambda_k\norm{F(x_k)-F(y_k)}
 \le\frac15e^{1/5}\norm{x_k-y_k}<\frac14\norm{x_k-y_k}.
\end{gather*}
Thus, Lemma~\ref{lem:energies} applies.
\item \emph{Change of the stepsize.} For $k\ge2$, the previous bounds give
\begin{gather*}
 \norm{F(y_k)}\le\frac54\cdot\frac65B_k=\frac32B_k,\\
 \norm{F(x_k)-F(y_k)}+\frac56\norm{F(y_k)}
 \le\Bigl(\frac3{10}+\frac54\Bigr)B_k=\frac{31}{20}B_k.
\end{gather*}
Let $t=\lambda_k/\Lambda_{k-1}$. Since $\lambda_k\le4\lambda_{k-1}\le4\Lambda_{k-1}$,
we have $t\le4$. By~\eqref{eq:maeg-identities} and the inequality
$\norm{u_{k-1}-y_{k-1}}\le2\Lambda_{k-1}B_k/3$,
\[
 \frac{B_{k+1}}{B_k}\le\frac32+\frac1{(1+3t/2)(1+t)}+\frac{93}{40}\cdot\frac t{1+t}<4,
\]
where for $t\le1/2$ the right-hand side is at most $3/2+1+31/40$, and for $1/2\le t\le4$
it is at most $3/2+8/21+93/50$.
For $k=1$,
\begin{gather*}
 \norm{F(y_1)}\le\frac65\norm{F(s)},\\
 \norm{u_1-y_1}=\lambda_1\Bigl\|\frac56F(s)-\frac{F(y_1)-F(s)}6\Bigr\|\le\frac{13}{15}\lambda_1\norm{F(s)},
\end{gather*}
hence $B_2\le\frac52\norm{F(s)}$. Both arguments of the minimum at the next iteration
are at least $\lambda_k/4$, that is, $\lambda_{k+1}\ge\lambda_k/4$.
\item \emph{Sum of the stepsizes.} The function $\log(1+t)/t$ is decreasing, hence
$\log(1+t)\ge(\log5/4)\,t\ge2t/5$ for $0\le t\le4$ and
\[
 \sum_{k=2}^K\frac{\lambda_k}{\Lambda_{k-1}}
 \le\frac52\log\frac{\Lambda_K}{\lambda_1}.
\]
Moreover,
$\lambda_k/\Lambda_{k-1}=(2/3+\lambda_k/\Lambda_{k-1})\alpha_k\le14\alpha_k/3$
and $\lambda_k\Lambda_{k-1}\le4\lambda_{k-1}\Lambda_{k-1}$.
By Lemma~\ref{lem:energies} and the Cauchy--Schwarz inequality,
\begin{align*}
 &\sum_{k=2}^K\lambda_k\norm{F(y_{k-1})}
 \le\Bigl(\sum_{k=2}^K\lambda_k\Lambda_{k-1}\norm{F(y_{k-1})}^2\Bigr)^{1/2}
 \Bigl(\sum_{k=2}^K\frac{\lambda_k}{\Lambda_{k-1}}\Bigr)^{1/2}
 \le6R\sqrt{\frac52\log\frac{\Lambda_K}{\lambda_1}},\\
 &\sum_{k=2}^K\frac{\lambda_k}{\Lambda_{k-1}}\norm{u_{k-1}-y_{k-1}}
 \le\Bigl(\frac{14}3\cdot\frac{37}3R^2\Bigr)^{1/2}
 \Bigl(\frac52\log\frac{\Lambda_K}{\lambda_1}\Bigr)^{1/2}
 \le12R\sqrt{\log\frac{\Lambda_K}{\lambda_1}}.
\end{align*}
Let $N_A$ be the number of iterations $k\ge2$ with
$\lambda_k=1/(6(L_0+L_1B_k))$. Then
\begin{align*}
 N_A&\le6\sum_{k=2}^K\lambda_k(L_0+L_1B_k)\\
 &\le6L_0\Lambda_K+6L_1R\Bigl(6\sqrt{\tfrac52}+\frac32\cdot12\Bigr)
 \sqrt{\log\frac{\Lambda_K}{\lambda_1}}
 \le6L_0\Lambda_K+165L_1R\sqrt{\log\frac{\Lambda_K}{\lambda_1}},
\end{align*}
since $6\sqrt{5/2}<19/2$.
On each of the remaining $K-1-N_A$ iterations, the stepsize increases fourfold,
and on iterations of the first type it decreases at most fourfold. Therefore,
\[
 (K-1-N_A)\log4-N_A\log4\le\log\frac{\lambda_K}{\lambda_1}
 \le\log\frac{\Lambda_K}{\lambda_1},
\]
whence
\[
 K\le1+2N_A+\frac1{\log4}\log\frac{\Lambda_K}{\lambda_1}.
\]
\item \emph{Termination.} As long as $\Lambda_{K-1}<3R/\varepsilon$, we have
\begin{gather*}
 \Lambda_K<\frac{3R}\varepsilon+\frac1{6L_0},\\
 \log\frac{\Lambda_K}{\lambda_1}\le\mathcal H.
\end{gather*}
Substitution gives
\[
 K\le1+2\Bigl(\frac{18L_0R}\varepsilon+1+165L_1R\sqrt{\mathcal H}\Bigr)+\frac{\mathcal H}{\log4}
 \le3+\frac{36L_0R}\varepsilon+330L_1R\sqrt{\mathcal H}+\mathcal H.
\]
Infinitely many iterations with $\Lambda_K<3R/\varepsilon$ would contradict
this bound. When $\Lambda_K\ge3R/\varepsilon$, Lemma~\ref{lem:energies}
guarantees termination.
\item \emph{Location of the point.} By the inequality
$(x+y)^2\le\frac65x^2+6y^2$ and Lemma~\ref{lem:energies},
\[
 \norm{y_K-z^\ast}^2
 \le\frac65\left(\norm{y_K-u_K}^2+5\norm{u_K-z^\ast}^2\right)
 \le\frac{74}5R^2<16R^2.\qedhere
\]
\end{itemize}
\end{proof}

\subsection{Final bound}

\begin{proof}[Proof of Theorem~\ref{thm:upper}\proofup{thm:upper}]
If \method{ABSEIL} outputs the point $s$ without running \method{BELAY}, then it
satisfies $\norm{F(s)}\le\varepsilon$ and $\norm{s-z^\ast}\le R$
by~\eqref{eq:stage-one}, and therefore $\norm{s-z_0}\le2R$. Otherwise, the output is the point
$y_K$ of the second stage, and by Proposition~\ref{prop:adaptive}
and~\eqref{eq:stage-one},
$\norm{y_K-z_0}\le\norm{y_K-z^\ast}+\norm{z^\ast-z_0}<4R+R=5R$.
We now bound the number of oracle calls. The first stage makes $J$ steps, and the second one makes
$K$ iterations, so the total number of calls is at most $1+J+2K$.
\begin{itemize}
\item \emph{Case $L_1R<1$.} Then $J=0$, $s=z_0$, and by
Lemma~\ref{lem:gronwall}(b),
$L_1\norm{F(z_0)}\le L_0(e^{L_1R}-1)\le2L_0$. Therefore,
$\mathcal H\le\log(3+54L_0R/\varepsilon)$ and
\[
 K\le3+\frac{36L_0R}\varepsilon+330\sqrt{\mathcal H}+\mathcal H
 \lesssim1+\frac{L_0R}\varepsilon.
\]
\item \emph{Case $L_1R\ge1$.} By Lemma~\ref{lem:long-step},
$J\le4(L_1R)^{6/5}$. By~\eqref{eq:stage-one},
$L_0+L_1\norm{F(s)}\le L_0\bigl(1+4L_1Re^{(L_1R)^{2/5}}\bigr)$, hence
\[
 \mathcal H\le\log\Bigl(1+\frac{18L_0R}\varepsilon\Bigr)+\log(1+4L_1R)+(L_1R)^{2/5}.
\]
Using $\sqrt{x+y+z}\le\sqrt x+\sqrt y+\sqrt z$ and the inequality
\[
 \log\Bigl(1+\frac{18L_0R}\varepsilon\Bigr)
 \le\log19+\log_+\frac{L_0R}\varepsilon
 \le\log19+\log(L_1R)+\frac{L_0R}{\varepsilon L_1R},
\]
which follows from $\log x\le x$ with $x=L_0R/(\varepsilon L_1R)$,
we obtain
\begin{align*}
 L_1R\sqrt{\mathcal H}
 &\le L_1R\sqrt{\log19+\log(L_1R)}+\sqrt{\frac{L_0R}\varepsilon L_1R}
 +L_1R\sqrt{\log(1+4L_1R)}+(L_1R)^{6/5}\\
 &\by{\lesssim}{(\mathrm a)}\frac{L_0R}\varepsilon+(L_1R)^{6/5},
\end{align*}
where (a) consists of the inequalities $\sqrt{xy}\le(x+y)/2$ and $\log(1+4x)\lesssim x^{2/5}$
for $x\ge1$. Similarly,
$\mathcal H\lesssim1+L_0R/\varepsilon+(L_1R)^{2/5}$. Therefore,
\[
 1+J+2K\lesssim1+\frac{L_0R}\varepsilon+(L_1R)^{6/5}.
\]
\end{itemize}
Finally, for $\varepsilon=\eta/(6R)$, we have $L_0R/\varepsilon=6L_0R^2/\eta$, and by
Lemma~\ref{lem:gap-norm},
\[
 \Gap_R(\hat z)\le\bigl(\norm{\hat z-z_0}+R\bigr)\norm{F(\hat z)}\le6R\cdot\frac\eta{6R}=\eta.\qedhere
\]
\end{proof}

\section{Proofs of the lower bounds}\label{app:lower}

The proofs follow the order of the construction in Section~\ref{sec:lower}: curvature and inversion, revelation of one coordinate at a time, the radius, smoothing and the test point, coordinate hiding, and then the linear term. Below, $f=Q+I_K$, $s=\one^\top v$ and $r=\norm v$.

\subsection{Curvature of Q and properties of the inverse operator}
\label{app:lower-chain}\label{app:lower-inverse}
Bounds~\eqref{eq:chain-bounds} and the inequalities $\vartheta\le A$, $\gamma v_i+(Bv)_i\le C$ from Section~\ref{sec:lower-face} are verified as follows.
\begin{itemize}
\item \emph{Gradient formula.} For $v\ne0$, the derivative of the first term of $Q$ with respect to each coordinate equals $2\log(1+\ell_1s/\ell_0)/\ell_1$. For the second term, $\partial_i r=v_i/r$ and $\partial_i s=1$ give
\[
\partial_i\left[r\log\frac{\sqrt n r}{s}\right]
=\frac{v_i}{r}\left(1+\log\frac{\sqrt n r}{s}\right)-\frac rs,
\]
which yields~\eqref{eq:perspective-gradient}.
\item \emph{Bounds on $\vartheta$ and $\gamma$.} On $K$ we have $r\le s\le G_1$ and $1\le\sqrt n r/s\le\sqrt n$. Therefore, the first term of $\vartheta$ is at most $2\log(1+\ell_1G_1/\ell_0)/\ell_1=A$, and the second one is nonpositive. Moreover, $0\le\gamma v_i\le\gamma r\le4\log(e\sqrt n)/\ell_1$, which implies~\eqref{eq:chain-bounds}.
\item \emph{Bound on $Bv$.} We have $|(Bv)_i|\le5\sqrt m\log(e\sqrt n)/\ell_1$.
The last bound follows from the fact that the sum of the coordinates of $v$ belonging to one
block does not exceed the weight of this block: $\sum_{\bl(k)=\bl(i)}v_k\le G_i$.
Together with the previous item, this gives $\gamma v_i+(Bv)_i\le C$.
\end{itemize}

\begin{proof}[Proof of Lemma~\ref{lem:curvature}\proofup{lem:curvature}]
\leavevmode
\begin{itemize}
\item \emph{Hessian of the second part.} Let $p=v/s$, $\alpha=\one^\top\xi$ and $\xi_\perp=\xi-p\alpha$. Then $\one^\top p=1$, $\norm p\le1$ and $\xi=p\alpha+\xi_\perp$. The second part of $Q$ equals $s\,a(p)$, where $a(p)=4\norm p\log(\sqrt n\norm p)/\ell_1$. For the radial function $a(p)=b(\norm p)$ with $r_p=\norm p$, we have
\[
\nabla^2a(p)=\frac{b'(r_p)}{r_p}\left(I-\frac{pp^\top}{r_p^2}\right)
+b''(r_p)\frac{pp^\top}{r_p^2}.
\]
Here $b'(r_p)=4(1+\log(\sqrt n r_p))/\ell_1$ and $b''(r_p)=4/(\ell_1r_p)$, hence
\[
\nabla^2a(p)=\frac4{\ell_1r_p}\left[I+\log(\sqrt n r_p)
\left(I-\frac{pp^\top}{r_p^2}\right)\right]\succeq\frac4{\ell_1r_p}I.
\]
The inequality holds because the logarithm is nonnegative and the matrix in parentheses after it is an orthogonal projector.
\item \emph{Second derivative along a line.} Let $v(t)=v+t\xi$, $s(t)=s+t\alpha$ and $p(t)=v(t)/s(t)$. Then $p'(0)=\xi_\perp/s$ and $p''(0)=-2\alpha\xi_\perp/s^2$. The second derivative of $s(t)a(p(t))$ at zero equals
\[
2\alpha\ip{\nabla a}{p'}+s\ip{\nabla^2a\,p'}{p'}+s\ip{\nabla a}{p''}
=\frac1s\ip{\nabla^2a\,\xi_\perp}{\xi_\perp},
\]
since the first and third terms cancel. Hence,
\[
D^2Q(v)[\xi,\xi]\ge
\frac{2\alpha^2}{\ell_0+\ell_1s}+\frac{4\norm{\xi_\perp}^2}{\ell_1\norm v}.
\]
\item \emph{Lower bound on the curvature.} Let $D_v=\ell_0+\ell_1s\norm p$. Since $0\le\norm p\le1$,
\[
D_v-\norm p^2(\ell_0+\ell_1s)
=\ell_0(1-\norm p^2)+\ell_1s(\norm p-\norm p^2)\ge0.
\]
Moreover, $4/(\ell_1\norm v)\ge2/D_v$. Therefore,
\[
D^2Q(v)[\xi,\xi]\ge
\frac{2\norm p^2\alpha^2+2\norm{\xi_\perp}^2}{D_v}
\by{\ge}{(\mathrm a)}\frac{\norm\xi^2}{\ell_0+\ell_1\norm v},
\]
where (a) is the inequality $\norm{a+b}^2\le2\norm a^2+2\norm b^2$.
\item \emph{Generalized strong monotonicity.} The notation $a\in\partial f(v)$ means $f(w)\ge f(v)+\ip a{w-v}$ for all $w$. The function $Q$ is convex by the Hessian bound just proved. It is nonnegative and continuous at zero: the first term is $O(s^2)$, and the second one is at most $4\norm v\log\sqrt n/\ell_1$. For $v,v'\in K$, set $V=\max\{\norm v,\norm{v'}\}$. On the segment between them the norm is at most $V$, so integrating the curvature gives the strong convexity inequality with coefficient $1/(\ell_0+\ell_1V)$. At a boundary point we use the one-sided derivative, and at zero the continuous extension along the segment. For any $a\in\partial f(v)$ and $a'\in\partial f(v')$, we obtain
\begin{align*}
f(v')&\ge f(v)+\ip a{v'-v}+\frac{\norm{v-v'}^2}{2(\ell_0+\ell_1V)},\\
f(v)&\ge f(v')+\ip{a'}{v-v'}+\frac{\norm{v-v'}^2}{2(\ell_0+\ell_1V)}.
\end{align*}
Adding these inequalities, we obtain
\begin{equation}\label{eq:inverse-strong}
\ip{a-a'}{v-v'}\ge\frac{\norm{v-v'}^2}{\ell_0+\ell_1V}.
\end{equation}
\item \emph{Existence of a solution.} Define
\[
\operatorname{prox}_f(x)=\arg\min_{w\in K}\{Q(w)+\tfrac12\norm{w-x}^2\}.
\]
The minimum exists by compactness and is unique because the quadratic term is strictly convex. The optimality condition $x-\operatorname{prox}_f(x)\in\partial f(\operatorname{prox}_f(x))$ and the monotonicity of the subdifferential imply that this map is nonexpansive and hence continuous. By Brouwer's theorem, the continuous map $v\mapsto\operatorname{prox}_f(v-Bv+u+\sigma)$ from $K$ into itself has a fixed point. This is exactly a solution of~\eqref{eq:inverse-relation}.
\item \emph{Uniqueness, monotonicity and Lipschitz continuity.} For solutions $v,v'$ corresponding to queries $u,u'$, substitute $a=u+\sigma-Bv$ and $a'=u'+\sigma-Bv'$ into~\eqref{eq:inverse-strong}. Skew-symmetry eliminates $B$:
\[
\ip{u-u'}{v-v'}\ge\frac{\norm{v-v'}^2}{\ell_0+\ell_1V}.
\]
For $u=u'$ we obtain $v=v'$, that is, uniqueness. For different queries this is the monotonicity of $\Phi$. By the Cauchy--Schwarz inequality, we also have
\[
\norm{\Phi(u)-\Phi(u')}\le
\bigl(\ell_0+\ell_1\max\{\norm{\Phi(u)},\norm{\Phi(u')}\}\bigr)\norm{u-u'}.
\]
The maximum at the endpoints does not exceed the supremum over the segment, so condition~\eqref{eq:smooth} holds. Finally, $\norm\Phi\le G_1$ and $\Phi(-\sigma)=0$, since $0\in\partial f(0)$.
\end{itemize}
\end{proof}

\subsubsection{KKT conditions from the geometry of the simplex}
At a point $v\ne0$ we have $\partial f(v)=\nabla Q(v)+N_K(v)$. Let $y\in N_K(v)$. The maximum of the linear function $\ip yw$ over $K$ equals
$\zeta=\max\{0,G_1y_1,\ldots,G_ny_n\}$, since it suffices to check the vertices. It is attained at $v$. Set $\lambda_i=\zeta/G_i-y_i\ge0$. Then
\[
\zeta=\ip yv=\zeta\sum_i v_i/G_i-\sum_i\lambda_i v_i,
\]
whence $\zeta(1-\sum_i v_i/G_i)+\sum_i\lambda_i v_i=0$. All terms are nonnegative, so each of them is zero. We obtain $y_i=\zeta/G_i-\lambda_i$, $\lambda_iv_i=0$ and the complementary slackness condition. Substituting into inclusion~\eqref{eq:inverse-relation} and taking~\eqref{eq:perspective-gradient} into account gives the KKT conditions: for some $\zeta\ge0$,
\begin{equation}\label{eq:kkt}
u_i+\sigma_i-\vartheta-\gamma v_i-(Bv)_i-\zeta/G_i
\begin{cases}=0,&v_i>0,\\\le0,&v_i=0,
\end{cases}
\end{equation}
and $\zeta(1-\sum_iv_i/G_i)=0$. If $\zeta>0$, the answer lies on the outer face.

\subsection{Revelation of one coordinate at a time}
\label{app:lower-reveal}
\begin{proof}[Proof of Lemma~\ref{lem:chain}\proofup{lem:chain}]
\leavevmode
\begin{itemize}
\item \emph{Part (a).} Let $|u_n|\le\delta$ and $v=\Phi(u)$. The answer $v=0$ is impossible. Otherwise, by the definition of the subgradient along $te_n$,
\[
A+2C-\delta\le u_n+\sigma_n\le\lim_{t\downarrow0}\frac{Q(te_n)}t
=\frac4{\ell_1}\log\sqrt n<C,
\]
which contradicts $A\ge0$ and $C>\delta$. If $\zeta=0$, then condition~\eqref{eq:kkt} for coordinate $n$ gives
\[
A+2C-\delta\le u_n+\sigma_n
\le\vartheta+\gamma v_n+(Bv)_n
\by{\le}{(\mathrm a)}A+C,
\]
where (a) is given by the bounds of Section~\ref{sec:lower-face}. This contradicts $C>\delta$. Hence, $\zeta>0$, the answer lies on the outer face, and by Lemma~\ref{lem:face} $\norm{\Phi(u)}\ge G_n/\sqrt n$.
\item \emph{Part (b), a jump within a block.} Let $|u_i|\le\delta$ for $i\ge j$, and let $l$ be the last nonzero coordinate of $v=\Phi(u)$. Suppose that $l>j$. By Part (a), the answer lies on the outer face. Assume first that the block of coordinate $l$ contains a coordinate $t\ge j$ preceding it, and let $\beta=5\sqrt m\log(e\sqrt n)/(\ell_1G_l)$. Subtract condition~\eqref{eq:kkt} for $t$ from the equality for $l$. The terms $\zeta/G_i$ cancel, since $G_t=G_l$. We obtain
\begin{align*}
-2\delta
&\by{\ge}{(\mathrm a)}(u_l+\sigma_l)-(u_t+\sigma_t)
\ge\gamma(v_l-v_t)+(Bv)_l-(Bv)_t\\
&\by{=}{(\mathrm b)}(\gamma+\beta)v_l+(\beta-\gamma)v_t+2\beta\sum_{t<i<l}v_i
\by{\ge}{(\mathrm c)}0,
\end{align*}
where (a) follows from the inequalities $\sigma_t-\sigma_l\ge4\delta$ and $|u_t|,|u_l|\le\delta$, (b) follows from the sign pattern of $B$, under which the terms before $t$ and after $l$ cancel, and (c) follows from the inequality $\gamma\le\beta$. The latter follows from~\eqref{eq:chain-bounds} and the bound
\[
1=\sum_{i\le l}\frac{v_i}{G_i}\le\norm v\Bigl(\sum_{i\le l}G_i^{-2}\Bigr)^{1/2}
\le\norm v\frac{\sqrt m}{G_l\sqrt{1-e^{-2}}}.
\]
This is impossible.
\item \emph{Part (b), a transition to another block.} Otherwise, $l$ is the first coordinate of a block strictly after the block containing $j$. After $l$ the answer is zero, so $(Bv)_l=0$. Let $r=G_l/G_j\le1/e$. Conditions~\eqref{eq:kkt} for $j$ and $l$ give $\zeta\ge G_j(A+C-\delta-\vartheta)$ and $\zeta\le G_l\bigl(A+2C+4(m-1)\delta+\delta-\vartheta\bigr)$. Comparing these bounds and dividing by $G_j$, we obtain
\begin{align*}
0&\ge(1-r)(A-\vartheta)+(1-2r)C-4r(m-1)\delta-(1+r)\delta\\
&\by{\ge}{(\mathrm a)}\frac1{\ell_1}\Bigl[\Bigl(1-\frac2e\Bigr)(5\sqrt m+4)-\frac{\sqrt m}e-\frac12\Bigr]
\by{>}{(\mathrm b)}0,
\end{align*}
where (a) follows from the inequalities $\vartheta\le A$, $r\le1/e$, $\log(e\sqrt n)\ge1$ and the definition of $\delta$, and (b) follows from the fact that the expression in brackets equals $(5-11/e)\sqrt m+7/2-8/e>0$. This is a contradiction.
\item \emph{Dependence on the prefix.} Hence, $l\le j$. The remaining coordinates of the answer solve the same inclusion on the first $j$ coordinates, where the solution is unique. Therefore, they depend only on $u_1,\ldots,u_j$. For $j=1$, only the first coordinate remains, and the equation defining the outer face gives $\Phi(u)=G_1e_1$.
\end{itemize}
\end{proof}

\subsection{Exact choice of the length and control of the radius}
\label{app:lower-radius}
\proofdown{lem:radius}
\begin{lemma}[chain length and radius]\label{lem:radius}
Let $h\ge2$ be an integer, $m=h^2$, $H=\lceil h\log(4h)\rceil$ and $n=mH$. Then
\begin{equation}\label{eq:chain-radius}
 \max\bigl\{\norm\sigma,\norm w\bigr\}\le\frac{128}{\ell_1}h^{5/2}\log^{3/2}(4h)
\end{equation}
and $h^3\log(4h)\le n\le2h^3\log(4h)$.
\end{lemma}

\begin{proof}[Proof of Lemma~\ref{lem:radius}\proofup{lem:radius}]
Let $\lambda=\log(4h)$.
\begin{itemize}
\item \emph{Chain length and logarithms.} We have $h^3\lambda\le n\le2h^3\lambda$, $\log(n+1)\le6\lambda$ and $\log(e\sqrt n)\le3\lambda$. The last bounds follow, for example, from $\lambda\le2h$ and $n\le4h^4$.
\item \emph{Bound on the shift.} The three parts of $\ell_1\sigma_i$ are bounded separately:
\begin{align*}
2\log(1+n^3e^{H-1})&\le6\log(n+1)+2H\le22h\lambda,\\
2(5h+4)\log(e\sqrt n)&\le44h\lambda,\\
4\ell_1(m-\pos(i))\delta&\le h\le h\lambda.
\end{align*}
Therefore, $\ell_1\max_i\sigma_i\le67h\lambda$.
\item \emph{Bound on the test point.} For the point $w$ from Lemma~\ref{lem:witness} and $i<n$, the gradient at the vertex is at least $-4/\ell_1$, and the term $G_n(Be_n)_i$ is either zero or equal to $-5h\log(e\sqrt n)/\ell_1$. All $w_i$ are negative. For $i=n$ this is seen from the explicit formula in the proof of Lemma~\ref{lem:witness}. Hence, $\ell_1\max_i|w_i|\le67h\lambda+4+15h\lambda\le83h\lambda$.
\item Multiplying by $\sqrt n\le\sqrt2\,h^{3/2}\lambda^{1/2}$ and using $83\sqrt2<128$, we obtain
\begin{equation}\label{eq:radius-detail}
\max\{\norm\sigma,\norm w\}\le\frac{128}{\ell_1}h^{5/2}\lambda^{3/2}.
\end{equation}
\end{itemize}
\end{proof}

\subsubsection{The largest admissible parameter and the resulting order}
Now $\ell_0=L_0/4$, $\ell_1=L_1/2$, $x=L_1R$, and $h\ge256$ is the largest integer such that $2048h^{5/2}\lambda^{3/2}\le x$.
\begin{itemize}
\item \emph{The solution and the test point lie in the ball.} From~\eqref{eq:radius-detail} we obtain $\norm\sigma\le R/8$ and $\norm w\le R/8$.
\item \emph{Lower bounds on $n$.} When passing from $h$ to $h+1$, the expression $h^{5/2}\log^{3/2}(4h)$ increases by a factor of less than $4$, so $x<8192h^{5/2}\lambda^{3/2}$. Hence,
\[
\frac{n^{5/2}}x>\frac{h^5\lambda}{8192}>4,
\]
\[
\frac{n^3}{x^2}>\frac{h^4}{8192^2}>4.
\]
Both strict bounds hold for $h\ge256$.
\item \emph{Order of $n$.} Also, $x^{6/5}=\Theta(h^3\lambda^{9/5})$ and $\log(2+x)=\Theta(\lambda)$, whence $n=\Theta\bigl(x^{6/5}/\log^{4/5}(2+x)\bigr)$.
\end{itemize}

\subsubsection{Verification of the coordinate hiding parameters}
For the norm we set $q=L_0R/\varepsilon$, and for the gap $q=L_0R^2/\eta$. We take $\mu=L_0/[64(1+q+x)]$. After smoothing, the constants are $a_0=L_0/2$ and $a_1=L_1/2$.
\begin{itemize}
\item \emph{Condition of Lemma~\ref{lem:localization}.} We have $\mu\le a_0/(8(1+a_1R))$.
\item \emph{Conditions of Lemma~\ref{lem:hardness}.} We have $\mu\le L_0/(64q)$, that is, $\mu\le\varepsilon/(16R)$ for the norm and $\mu\le\eta/(16R^2)$ for the gap. For the norm, $M=G_1\ge G_n/\sqrt n>2\varepsilon$. For the gap, $G_n/\ell_1>4\eta$, so in the regime of large $x$ under consideration, $M\ge G_n>2L_1\eta\ge\eta/R$.
\item \emph{Solution.} The solution after smoothing has norm at most $R/8+R/16<R/4$.
\end{itemize}

\subsection{Smoothing and the test point}
\label{app:lower-smoothing}

\begin{proof}[Proof of Lemma~\ref{lem:smoothing}\proofup{lem:smoothing}]
\leavevmode
\begin{itemize}
\item \emph{Monotonicity and boundedness.} The convolution $\Phi_\rho$ is smooth. As an average of vectors from $K$, it lies in $K$, so its norm is at most $G_1$. Monotonicity follows from
\[
\ip{\Phi_\rho(u)-\Phi_\rho(u')}{u-u'}
=\int\ip{\Phi(u-y)-\Phi(u'-y)}{u-u'}k_\rho(y)dy\ge0.
\]
\item \emph{Solution.} If $u\le-\sigma$ coordinatewise, then $\ip{u+\sigma}v\le0\le Q(v)$ on $K$, so $u+\sigma\in\partial f(0)$ and $\Phi(u)=0$. For $u=-\sigma-\rho\one$, all $u-y$ on the support of the kernel satisfy this condition. Therefore, $\Phi_\rho(u)=0$.
\item \emph{Chain properties.} If $|u_i|\le\theta=\delta/2$ for $i\ge j$, then $|u_i-y_i|\le\theta+\rho\le\delta$. Each answer under the integral has a zero tail after $j$, depends only on the prefix, and lies on the outer face when $u_n$ is small. Integration preserves all three properties.
\item \emph{$(L_0,L_1)$-Lipschitz continuity.} Set $K_{\mathrm{Lip}}=\ell_0+\ell_1G_1$. The original $\Phi$ is globally $K_{\mathrm{Lip}}$-Lipschitz, and almost everywhere $\norm{D\Phi}\le\ell_0+\ell_1\norm\Phi$ by Lemma~\ref{lem:class}. For $y$ in the support,
\[
\norm{\Phi(u-y)-\Phi_\rho(u)}
\le\int K_{\mathrm{Lip}}\norm{y-t}\,k_\rho(t)dt\le2K_{\mathrm{Lip}}\rho.
\]
Differentiating the convolution gives
\[
\norm{D\Phi_\rho(u)}
\le\ell_0+\ell_1\norm{\Phi_\rho(u)}+2\ell_1K_{\mathrm{Lip}}\rho
\le2\ell_0+\ell_1\norm{\Phi_\rho(u)}.
\]
\end{itemize}
\end{proof}

\begin{proof}[Proof of Lemma~\ref{lem:witness}\proofup{lem:witness}]
\leavevmode
\begin{itemize}
\item \emph{The answer in a neighborhood of $w$.} At the vertex $G_ne_n$, the normal cone consists of the $y$ for which $y_n\ge0$ and $y_iG_i\le y_nG_n$ for $i<n$. This follows from comparing the values of the linear function at all vertices. Let $r_w=G_n/[2\ell_1(G_1+G_n)]$ and $\norm\xi\le r_w$. Then $y=e_n/\ell_1+\xi$ belongs to this cone, since $y_n\ge1/\ell_1-r_w>0$ and $y_iG_i\le r_wG_1\le(1/\ell_1-r_w)G_n\le y_nG_n$.
Consequently, $\Phi(w+\xi)=G_ne_n$ in the whole ball. Since $\rho\le r_w$, we obtain $\Phi_\rho(w)=G_ne_n$.
\item \emph{The last coordinate of $w$.} We have
\[
\ell_1w_n=2\log\frac{1+n^3}{1+n^3e^{H-1}}+4\log\sqrt n+1
-2(5\sqrt m+4)\log(e\sqrt n)\le-3,
\]
since the first logarithm is nonpositive and $4\log\sqrt n<4\log(e\sqrt n)$.
Therefore, for $|u_n|\le\delta$, we have $\ip{\Phi_\rho(w)}{u-w}=G_n(u_n-w_n)\ge G_n(3/\ell_1-\delta)\ge G_n/\ell_1$.
\item The bound $\norm w\le R/8$ for the parameters of Section~\ref{sec:lower-assemble} is obtained in Appendix~\ref{app:lower-radius}.
\end{itemize}
\end{proof}

\subsection{Coordinate hiding preserves the class and the success condition}
\label{app:lower-localization}

\begin{proof}[Proof of Lemma~\ref{lem:localization}\proofup{lem:localization}]
Denote the original $(L_0,L_1)$-Lipschitz constants by $a_0,a_1$. Since $F$ is bounded, condition~\eqref{eq:smooth} makes it globally Lipschitz with constant $a_0+a_1M$. The derivatives below exist almost everywhere.
\begin{itemize}
\item \emph{Strong monotonicity.} We have $D\widetilde F=\chi DF+F\nabla\chi^\top+\mu I$. The monotonicity of $F$ gives $\ip{DF\xi}\xi\ge0$ almost everywhere. The middle term has norm at most $2M/L$, hence
\[
\ip{D\widetilde F\xi}\xi\ge(\mu-2M/L)\norm\xi^2\ge\tfrac\mu2\norm\xi^2.
\]
Integration over almost all parallel segments and passage to the limit, as in Lemma~\ref{lem:class}, give $\mu/2$-strong monotonicity.
\item \emph{Solution.} A zero $z_\mu^*$ of the operator $F+\mu I$ exists: the map $z\mapsto-F(z)/\mu$ maps the ball of radius $M/\mu$ into itself and has a fixed point by Brouwer's theorem. It is unique by strong monotonicity. By the monotonicity of the original $F$,
\[
0\le\ip{F(z_\mu^*)-F(z^*)}{z_\mu^*-z^*}
=-\mu\ip{z_\mu^*}{z_\mu^*-z^*}.
\]
Hence, $\norm{z_\mu^*}\le\norm{z^*}\le R<L$, so $\chi=1$ at this point, and it is also a zero of $\widetilde F$.
\item For $\norm z\le2R$ we have $\widetilde F=F+\mu z$ and
\[
\norm{D\widetilde F}\le a_0+a_1\norm{\widetilde F}+\mu(1+2a_1R)
\le\tfrac54a_0+a_1\norm{\widetilde F}.
\]
\item For $2R\le\norm z\le L$, monotonicity with respect to $z^*$ gives $\ip{F(z)}z\ge\ip{F(z)}{z^*}\ge-R\norm{F(z)}$. For $a=\norm{F(z)}$ and $b=\mu\norm z$ we obtain $2\mu Ra\le ab$ and
\begin{equation}\label{eq:noncancellation}
\norm{F(z)+\mu z}^2\ge a^2+b^2-ab\ge(a+b)^2/4.
\end{equation}
The last inequality is equivalent to $3(a-b)^2\ge0$. Therefore, $\norm F\le2\norm{\widetilde F}$ and $\norm{D\widetilde F}\le a_0+\mu+2a_1\norm{\widetilde F}$.
\item For $L\le\norm z\le2L$ we have $\norm{\widetilde F}\ge\mu\norm z-M\ge31M$ and
\[
\norm{D\widetilde F}\le a_0+a_1M+\frac{2M}L+\mu\le a_0+a_1M+\frac{3\mu}2.
\]
\item For $\norm z\ge2L$ we have $D\widetilde F=\mu I$.
\item In all four cases, $\norm{D\widetilde F}\le2a_0+2a_1\norm{\widetilde F}$ almost everywhere. By Lemma~\ref{lem:class}, this implies condition~\eqref{eq:smooth} with constants $2a_0,2a_1$. If the original $F\in C^2$, then also $\widetilde F\in C^2$: in a neighborhood of zero the factor $\chi$ is constant, so the nondifferentiability of the norm at zero causes no difficulty.
\end{itemize}
\end{proof}

\begin{proof}[Proof of Lemma~\ref{lem:hardness}\proofup{lem:hardness}]
\leavevmode
\begin{itemize}
\item \emph{Part (a).} Let $\mu\le\varepsilon/(16R)$, $M\ge\varepsilon$ and $\norm{\widetilde F(z)}\le\varepsilon$. For $L\le\norm z\le2L$ the norm of $\widetilde F$ is at least $31M>\varepsilon$, and for $\norm z\ge2L$ it is at least $64M>\varepsilon$. Therefore, $\norm z<L$. If $\norm z\le2R$, then $\norm{F(z)}\le\varepsilon+\mu\norm z\le\varepsilon+\varepsilon/8<2\varepsilon$. If $2R\le\norm z<L$, then by~\eqref{eq:noncancellation} $\norm{F(z)}\le2\norm{\widetilde F(z)}\le2\varepsilon$.
\item \emph{Part (b), bound on $\norm z$.} Let $\norm{z^*}\le R/4$, $\mu\le\eta/(16R^2)$, $M\ge\eta/R$ and $\Gap_R^{\widetilde F}(z)\le\eta$. Let $r=\norm{z-z_\mu^*}\ge R/2$ and $u=z_\mu^*+(R/2)(z-z_\mu^*)/r$. Then $\norm u\le3R/4$, and $\widetilde F=F+\mu I$ on the relevant segment. Therefore,
\begin{align*}
\eta&\ge\ip{\widetilde F(u)}{z-u}
=\frac{r-R/2}{R/2}\ip{\widetilde F(u)-\widetilde F(z_\mu^*)}{u-z_\mu^*}\\
&\ge\frac{r-R/2}{R/2}\,\mu R^2/4
=\frac{\mu R}{2}(r-R/2).
\end{align*}
Hence, $r\le R/2+2\eta/(\mu R)$, and for $r<R/2$ this bound also holds. Consequently,
\begin{equation}\label{eq:output-localization}
\norm z\le3R/4+2\eta/(\mu R)
\le3L/16+2M/\mu\le L/4.
\end{equation}
\item \emph{Part (b), bound on $\Gap^F_R(z)$.} In the ball of radius $R$ we have $\widetilde F(u)=F(u)+\mu u$, so
\begin{align*}
\Gap_R^F(z)&\le\Gap_R^{\widetilde F}(z)+\mu R\norm z+\mu R^2\\
&\le\eta+\mu R\left(3R/4+2\eta/(\mu R)\right)+\mu R^2\\
&=3\eta+\tfrac74\mu R^2<4\eta.
\end{align*}
\end{itemize}
\end{proof}

\subsection{Random embedding}
\label{app:lower-hiding}
\begin{proof}[Proof of Proposition~\ref{prop:random-hiding}\proofup{prop:random-hiding}]
\leavevmode
\begin{itemize}
\item \emph{Membership in the class.} Denote $F_0(z)=U\Phi(U^\top z)$. Since $U^\top U=I$, we have $\norm{F_0(z)}=\norm{\Phi(U^\top z)}$, $\norm{U^\top(z-z')}\le\norm{z-z'}$ and
\[
\ip{F_0(z)-F_0(z')}{z-z'}=\ip{\Phi(U^\top z)-\Phi(U^\top z')}{U^\top z-U^\top z'}\ge0.
\]
Therefore, $F_0$ is monotone, satisfies condition~\eqref{eq:smooth} with the same constants, is bounded by $M$, and has the zero $Uu^*$ with norm $\norm{u^*}\le R$. By Lemma~\ref{lem:localization}, the operator $F_U$ belongs to $\mathcal P_d$ with constants $2L_0$, $2L_1$.
\item \emph{The gap under the embedding.} Restricting the points in the supremum to vectors of the form $Uu$, we obtain
\[
\Gap_R^{F_0}(z)\ge\sup_{\norm u\le R}\ip{U\Phi(u)}{z-Uu}
=\Gap_R^\Phi(U^\top z).
\]
\item \emph{The ideal process.} Fix the randomness of the algorithm. After $t$ revelations, the auxiliary oracle knows only the first $t$ positions. For a query with $\norm z\ge2L$, it returns $\mu z$ without revealing coordinates. For a query with $\norm z<2L$, it reveals $\pi(t+1)$ and returns $\chi(\norm z/L)U\Phi(\widetilde u)+\mu z$, where
\[
\widetilde u_i=
\begin{cases}
z_{\pi(i)},&i\le t,\\
\operatorname{clip}_{[-\theta,\theta]}(z_{\pi(t+1)}),&i=t+1,\\
0,&i>t+1,
\end{cases}
\]
and $\operatorname{clip}_{[-\theta,\theta]}(a)=\max\{-\theta,\min\{a,\theta\}\}$. By condition~(a), the vector $\Phi(\widetilde u)$ has a zero tail after $t+1$, so the ideal answer depends only on the revealed positions. Clipping is needed for this to hold on any history, including histories with a large value of the new coordinate. If, on the other hand, all unrevealed coordinates in the query are at most $\theta$ in absolute value, then by condition~(a) the ideal answer coincides with the real one.
\item \emph{Probability of guessing.} A vector $z$ with $\norm z<2L$ has at most $4L^2/\theta^2$ coordinates exceeding $\theta$ in absolute value: each such coordinate contributes more than $\theta^2$ to $\norm z^2<4L^2$. In the ideal process, conditionally on the past history, the unrevealed positions are uniform among the at least $d-n$ remaining ones. Therefore, for one unrevealed position,
\[
\Pr\{|z_{\pi(i)}|>\theta\mid\text{ideal history}\}
\le\frac{4L^2}{(d-n)\theta^2}.
\]
We call guessing the following event: among the first $\min\{\tau,N\}$ queries of the ideal process and, if $\tau\le N$, the returned point, there is a point $z$ with $\norm z<2L$ and an unrevealed position $\pi(i)$ such that $|z_{\pi(i)}|>\theta$. The algorithm may not have queried the returned point, but this point can also guess a hidden coordinate. Summing over at most $n$ unrevealed positions and at most $N+1$ points, we obtain
\[
p_{\rm bad}\le\frac{4(N+1)nL^2}{(d-n)\theta^2}\le\frac1{12}.
\]
We condition on the ideal history, without the additional condition ``there has been no guessing before''. Otherwise, the uniformity of the unrevealed positions could be violated.
\item \emph{Success is impossible without guessing.} Suppose there was no guessing and $\tau\le N$. Then the real answers coincide with the ideal ones, and at most $N<n$ positions are revealed, so the position $\pi(n)$ is not revealed. If $\norm{\hat z}<2L$, then $|\hat z_{\pi(n)}|\le\theta$, and by condition~(b) $\norm{U\Phi(U^\top\hat z)}>2\varepsilon$. If success means a small gap, then the previous item gives $\Gap_R^{F_0}(\hat z)\ge\Gap_R^\Phi(U^\top\hat z)>4\eta$. In both cases, by Lemma~\ref{lem:hardness}, the point $\hat z$ is not a successful output. If $\norm{\hat z}\ge2L$, the same lemma rules out success. Therefore, $\Pr\{\text{success}\}\le\Pr\{\tau>N\}+p_{\rm bad}$.
\item \emph{Number of queries.} Let $T=\sup_U\E\tau(F_U)$. After averaging over $U$ and the randomness of the algorithm, Markov's inequality gives
\[
\frac23\le\Pr\{\tau>N\}+p_{\rm bad}
\le\frac{T}{N+1}+\frac1{12},
\]
whence $T\ge7(N+1)/12$. The distribution of the embedding is chosen independently of the randomness of the algorithm, so there exists a fixed matrix $U$ on which the algorithm makes at least $7(N+1)/12$ oracle calls in expectation.
\end{itemize}
\end{proof}

\subsection{Linear systems: accuracy, probability and expected running time}
\label{app:lower-linear}
\begin{proof}[Proof of Lemma~\ref{lem:linear-budget}\proofup{lem:linear-budget}]
\leavevmode
\begin{itemize}
\item \emph{External result.} For $\kappa\ge2$ and a sufficiently large dimension $d$, achieving relative residual $1/d$ with probability $5/6$ for all symmetric matrices with spectrum in $[-\kappa,-1]\cup[1,\kappa]$ requires $\Omega(\kappa\log d)$ arbitrary adaptive multiplications under a hard budget. This is the form of Theorems 2.6--2.8 of \citet{derezinski2026} that we use. The spectral subclass is preserved in their proof.
\item \emph{Guaranteed success.} Suppose that a solver for relative residual $1/2$ succeeds with probability $2/3$ and uses at most $T$ multiplications in expectation on any system. After each run, we check the residual with an additional multiplication and repeat an independent run in case of failure. The probability of reaching the $j$th run is at most $(1/3)^{j-1}$, so the expected cost of a verified success is at most
\[
\sum_{j\ge1}(1/3)^{j-1}(T+1)=\tfrac32(T+1).
\]
The cost of a run may depend on its success. What matters is only that the event of reaching it depends on the previous runs.
\item \emph{Refinement.} Start with $x=0$ and $r=b$. We solve $Au=r$ with the verified procedure and update $x\leftarrow x+u$, $r\leftarrow r-Au$. The residual norm decreases by at least a factor of two. After $\lceil\log_2d\rceil$ corrections, the relative residual is at most $1/d$, and the expected cost is at most $C(T+1)\log d$. If $r=0$ earlier, we stop.
\item \emph{Hard budget.} We truncate the whole procedure after $\lceil6C(T+1)\log d\rceil$ multiplications. By Markov's inequality, the probability of not finishing in time is at most $1/6$. The external result gives $C'(T+1)\log d\ge c\kappa\log d$, whence $1+T\gtrsim\kappa$. For $1\le\kappa<2$, the bound follows from $1+T\ge1$. The dimension can be chosen sufficiently large, since the complexity is taken uniformly over all dimensions.
\end{itemize}
\end{proof}

\begin{proof}[Proof of Proposition~\ref{prop:linear-lower}\proofup{prop:linear-lower}]
\leavevmode
\begin{itemize}
\item \emph{Membership in the class.} For the reductions of Section~\ref{sec:linear-lower}, the linear part of the operator is $S=a\begin{psmallmatrix}0&A\\-A&0\end{psmallmatrix}$. Since $A=A^\top$, we have $S^\top=-S$ and $\norm S=a\norm A$. This verifies monotonicity and $(L_0,L_1)$-Lipschitz continuity. Two queries to $A$ simulate one value of $F$.
\item \emph{Bound for the norm.} For $a=2\varepsilon/R$ and $s=R$, $\norm{F(x,y)}\le\varepsilon$ implies $\norm{A(\norm b/R)x-b}\le\norm b/2$. Therefore, Lemma~\ref{lem:linear-budget} gives the first bound of the proposition.
\item \emph{Gap of an affine operator.} For $F(z)=Sz+c$ with $S^\top=-S$, we have $\ip{F(u)}{z-u}=\ip c z+\ip u{-Sz-c}$. The maximum over the ball equals $\ip c z+R\norm{F(z)}$. Since $\ip{Sz}z=0$ and $F(z)=S(z-z^*)$, we have $\ip c z=\ip{F(z)}z=\ip{F(z)}{z^*}$. If $\norm{z^*}\le R/4$, we obtain $\Gap_R^F(z)\ge3R\norm{F(z)}/4$.
\item \emph{Bound for the gap.} For $a=16\eta/R^2$ and $s=R/4$, a successful output therefore gives
\[
\left\|A\frac{4\norm b}{R}x-b\right\|
\le\frac{4\norm b}{R}\frac{\norm{F(x,y)}}a
\le\frac{4\norm b}{R}\frac{4\eta}{3R}\frac{R^2}{16\eta}
=\frac{\norm b}{3}.
\]
For the norm and the gap, we apply Lemma~\ref{lem:linear-budget} to the resulting solver, which uses at most $2T$ multiplications in expectation. When the chosen $\kappa<1$, the corresponding ratio satisfies $L_0R/\varepsilon<2$ or $L_0R^2/\eta<16$, so the constant term $1$ suffices for the corresponding lower bound.
\end{itemize}
\end{proof}

\section{Bounded set}\label{app:constraints}

\subsection{Relation between the criteria}

Projection methods control the residual
$r_\gamma(z)=\norm{z-P_Z(z-\gamma F(z))}/\gamma$, while the natural measure of
nonoptimality is the tangential residual
$r_{\rm tan}(z)=\dist\bigl(0,F(z)+N_Z(z)\bigr)$.

\begin{lemma}[gaps and residuals]\label{lem:residuals}
Let $Z$ be a convex compact set of diameter $D$, $z\in Z$, $\gamma>0$, and
$p=P_Z(z-\gamma F(z))$. Then
\begin{gather*}
 r_\gamma(z)\le r_{\rm tan}(z),\\
 \Gap^{\rm S}_Z(z)\le D\,r_{\rm tan}(z),\\
 \gamma\,r_\gamma(z)^2\le\Gap^{\rm S}_Z(z).
\end{gather*}
If $F$ is Lipschitz on $[z,p]$ with constant $L_{zp}$, then
$\xi_p=(z-p)/\gamma-F(z)\in N_Z(p)$ and
$\norm{F(p)+\xi_p}\le(1+\gamma L_{zp})\,r_\gamma(z)$.
\end{lemma}

\begin{proof}
Let $\xi\in N_Z(z)$. Then $P_Z(z+\gamma\xi)=z$, and by nonexpansiveness of
the projection,
\[
 \gamma r_\gamma(z)=\norm{P_Z(z+\gamma\xi)-P_Z(z-\gamma F(z))}\le\gamma\norm{F(z)+\xi}.
\]
For $u\in Z$, we have $\ip\xi{z-u}\ge0$, hence
$\ip{F(z)}{z-u}\le\ip{F(z)+\xi}{z-u}\le D\norm{F(z)+\xi}$. Minimizing over $\xi$,
we obtain the first two inequalities. The projection inequality at the point $u=z$ gives
$\ip{z-\gamma F(z)-p}{z-p}\le0$, that is,
$\gamma\Gap^{\rm S}_Z(z)\ge\gamma\ip{F(z)}{z-p}\ge\norm{z-p}^2$.
Finally, $\xi_p\in N_Z(p)$ by the definition of the projection, and
\[
 \norm{F(p)+\xi_p}\le\norm{F(p)-F(z)}+\frac{\norm{z-p}}\gamma
 \le(L_{zp}+1/\gamma)\norm{z-p}.\qedhere
\]
\end{proof}

The reverse inequality $r_{\rm tan}(z)\lesssim r_\gamma(z)$ fails: for
$Z=[0,1]$, $F\equiv1$, and $0<z<\gamma$, we have $r_{\rm tan}(z)=1$ and
$r_\gamma(z)=z/\gamma$. The last claim of the lemma concerns the point $p$,
not $z$.

\subsection{Verification of the example on a bounded set}
\label{app:constraints-lower}
\begin{proof}[Proof of Theorem~\ref{thm:constrained-impossible}\proofup{thm:constrained-impossible}]
Let $\Gamma\ge72\eta/D^2$, $\kappa=\Gamma D^2/(72\eta)$, and $A=A^\top$ with spectrum in $[-\kappa,-1]\cup[1,\kappa]$. Denote $a=72\eta/D^2$, $s=Db/(12\norm b)$, and $x^*=A^{-1}s$.
The set $Z=B(0,D/3)\times[0,D/3]\times B(0,D/3)$ has diameter
$\sqrt{(2D/3)^2+(D/3)^2+(2D/3)^2}=D$.

The operator $F(x,t,y)=(aAy,\Gamma/L_1,a(s-Ax))$ is monotone: its linear part is skew-symmetric. On the whole space, $\norm{DF}=a\norm A\le\Gamma\le L_1\norm F\le L_0+L_1\norm F$.
The norm of $x^*$ is at most $D/12$. For the point $z^*=(x^*,0,0)$, we obtain $F(z^*)=(0,\Gamma/L_1,0)$ and $\ip{F(z^*)}{z-z^*}=\Gamma t/L_1\ge0$ for all $z\in Z$.

For a skew-symmetric linear part, $\ip{F(u)}{z-u}=\ip{F(z)}{z-u}$. Hence, the Minty gap coincides with the strong gap. Over the product of balls and a segment, the supremum is computed coordinatewise:
\begin{align*}
\Gap_Z^{\rm M}(x,t,y)
&=a\ip{Ax^*}y+\Gamma t/L_1
+\frac{aD}{3}\norm{Ay}+\frac{aD}{3}\norm{Ax-s}\\
&\ge\frac{aD}{3}\norm{Ax-s}=\frac{24\eta}{D}\norm{Ax-s}.
\end{align*}
The last line uses $t\ge0$ and
$\ip{Ax^*}y=\ip{x^*}{Ay}\ge-(D/12)\norm{Ay}$: the negative term is absorbed by the positive term $D\norm{Ay}/3$.
A gap $\le\eta$ implies $\norm{Ax-s}\le D/24$, and therefore
\[
\left\|A\frac{12\norm b}{D}x-b\right\|
\le\frac{12\norm b}{D}\frac D{24}=\frac{\norm b}{2}.
\]
Two multiplications by $A$ simulate one query. Therefore, $1+2T\gtrsim\kappa$.

To remove the additive one for bounded $\kappa$, we take $A=I$ and $b=\pm e_1$ in the same dimension. A successful answer requires $\norm{x\mp(D/12)e_1}\le D/24$, respectively. These two balls are disjoint. Before the first query, the algorithm cannot distinguish the instances. If the probability of making at least one query equals $p$, the average success probability over the two instances is at most $p+(1-p)/2$. Success probability $\ge2/3$ requires $p\ge1/3$, hence $T\ge1/3$. Together with the large-$\kappa$ bound, this gives
$T\gtrsim\kappa=\Theta(\Gamma D^2/\eta)$ for all $\kappa\ge1$.

The exact projection onto $Z$ does not help to extract information about $A,b$: the geometry of $Z$ does not depend on them. Since $\Gamma$ is arbitrarily large, the complexity of the class without control of the norm at the solution is infinite.

\end{proof}

\clearpage
\section{Summary of known bounds}\label{app:literature}
\relpenalty=500 \binoppenalty=700

\begin{center}
\footnotesize
\captionof{table}{Monotone operators: number of calls to the operator oracle.}
\label{tab:literature}
\renewcommand\arraystretch{1.1}\setlength\tabcolsep{4pt}
\newcolumntype{M}[1]{>{\raggedright\arraybackslash}m{#1}}
\renewcommand\tabularxcolumn[1]{m{#1}}
\begin{tabularx}{\linewidth}{@{}M{.27\linewidth}M{.13\linewidth}M{.19\linewidth}Y@{}}
\toprule
Method and paper & Smoothness & Criterion & Bound\\
\midrule
\multicolumn{4}{@{}l}{\textit{Lipschitz operator}}\\
\method{Mirror-Prox}\textsuperscript{\dag}\newline\citet[Proposition~2.2]{nemirovski2004}
 & $L$ & Minty gap of the averaged iterate on~$Z$
 & $\displaystyle O\Bigl(1+\frac{LD^2}{\eta}\Bigr)$\\
Lower bound\textsuperscript{\ddag}\newline\citet{nemirovsky1992}
 & $L$, linear & Minty gap on a ball
 & $\displaystyle \Omega\Bigl(\frac{LD^2}{\eta}\Bigr)$\\
Halpern iteration\newline\citet[Theorem~2.10]{diakonikolas2020}
 & $L$ & $\norm{F(\hat z)}\le\varepsilon$
 & $\displaystyle O\Bigl(\frac{LR}{\varepsilon}\log\frac{LR}{\varepsilon}\Bigr)$\\
\method{EAG}\textsuperscript{\S}\newline\citet[Theorems~1--3]{yoon2021}
 & $L$ & $\norm{F(z_k)}\le\varepsilon$
 & $\displaystyle O\Bigl(1+\frac{LR}{\varepsilon}\Bigr)$\\
\midrule
\multicolumn{4}{@{}l}{\textit{$(L_0,L_1)$-Lipschitz operator under $p$-quasi-sharpness}}\\
\method{EG}, Popov's method\textsuperscript{\P}\newline\citet[Theorems~5.2, 5.3, 6.3]{vankov2024},\newline
\method{PF-AEG}\textsuperscript{\P}\newline\citet[Theorems~3.4,~3.5]{chu2026}
 & $\alpha$-sym.\ & $\dist^2(z_k,U^\ast)$
 & linear rate for $p\le2$, $O\bigl(k^{-2/(p-2)}\bigr)$ for $p>2$\\
\midrule
\multicolumn{4}{@{}l}{\textit{$(L_0,L_1)$-Lipschitz operator without additional assumptions}}\\
\method{EG}\newline\citet[Theorem~3.6]{choudhury2025},\newline
\method{PolyakEG}\newline\citet[Theorem~3.4]{yoon2026polyak},\newline
\method{EG}, \method{PEG}\textsuperscript{\dag\dag}\newline\citet[Theorems~3.1, 3.2]{greenberg2026}
 & sym.\ $(L_0,L_1)$ & $\min_{k\le K}\norm{F(x_k)}\le\varepsilon$
 & $\displaystyle O\Bigl(1+\Bigl(\frac{L_0R}{\varepsilon}+L_1R\Bigr)^2\Bigr)$\\
\method{MAEG}\textsuperscript{\ddag\ddag}\newline\citet{sun2026}
 & continuous & $\norm{F(y_k)}\le\varepsilon$
 & no bound in terms of $L_0,L_1,R,\varepsilon$\\
LongTails\textsuperscript{\S\S}\newline\citet{tyurin2026longtails}
 & sym.\ $(L_0,L_1)$ & $\norm{F(\hat z)}\le\varepsilon$
 & $\displaystyle O\Bigl(1+\frac{L_0R}{\varepsilon}+(L_1R)^{4/3}\Bigr)$\\
\midrule
\method{ABSEIL} and lower bound\textsuperscript{\P\P}\newline\textbf{This paper}, Theorem~\ref{thm:main}
 & sym.\ $(L_0,L_1)$ & $\norm{F(\hat z)}\le\varepsilon$
 & $\displaystyle O\Bigl(1+\frac{L_0R}{\varepsilon}+(L_1R)^{6/5}\Bigr)$,\newline
 $\displaystyle\Omega\Bigl(1+\frac{L_0R}{\varepsilon}+\frac{(L_1R)^{6/5}}{\log^{4/5}(2+L_1R)}\Bigr)$\\
\bottomrule
\end{tabularx}
\par\smallskip\raggedright\scriptsize
\emph{Notation.} $R=\norm{z_0-z^\ast}$, $D$ is the diameter of the compact set $Z$,
$U^\ast$ is the solution set, and $O(\cdot)$ and $\Omega(\cdot)$ hide absolute constants.
The operator is monotone and the problem is posed on the whole $\R^d$ unless stated otherwise.\par
\emph{Smoothness.} $L$: the operator is $L$-Lipschitz.
sym.: symmetric $(L_0,L_1)$-Lipschitz continuity~\eqref{eq:smooth},
$\norm{F(u)-F(v)}\le\bigl(L_0+L_1\sup_{w\in[u,v]}\norm{F(w)}\bigr)\norm{u-v}$.
For $C^1$ operators, this is equivalent to $\norm{DF(z)}\le L_0+L_1\norm{F(z)}$.
$\alpha$-sym.: $\alpha$-symmetric generalized smoothness \citep{vankov2024},
which coincides with sym.\ for $\alpha=1$.\par\smallskip
\textsuperscript{\dag}~On a compact set with exact projections. The method uses as many projections as oracle calls.\quad
\textsuperscript{\S}~A lower bound of the same order holds for methods satisfying the linear-span assumption.\par
\textsuperscript{\ddag}~For deterministic methods, as stated by \citet[p.~242]{nemirovski2004}:
for a budget $N\le d/4$, there are bilinear problems with gap $\Omega(LD^2/N)$ on balls of radius $D$.\par
\textsuperscript{\P}~Under the $p$-quasi-sharpness condition $\ip{F(u)}{u-u^\ast}\ge\mu\dist(u,U^\ast)^p$, $\mu>0$.
It does not follow from monotonicity: for a skew-symmetric operator, the left-hand side is zero.
In \citet{chu2026}, the constants depend on the trajectory.\par
\textsuperscript{\dag\dag}~Constraints are allowed. For the Minty gap of the averaged iterate, the stepsize depends on $\norm{F(z^\ast)}$
and on the factor $e^{L_1\norm{z_0-z^\ast}}$.\par
\textsuperscript{\ddag\ddag}~Trial stepsizes are checked locally. For accepted stepsizes,
$\norm{F(y_k)}\le2R/\bigl((1-2\rho)\Lambda_k\bigr)$, but $\Lambda_k$ is not bounded in terms of $L_0,L_1,R,\varepsilon$.\par
\textsuperscript{\S\S}~Theorem \texttt{RadiusGeneralizedMonotone}.\quad
\textsuperscript{\P\P}~The upper bound is deterministic, and the lower bound holds for randomized
methods when $L_0R/\varepsilon\ge1$. The same holds for $\Gap_R(\hat z)\le\eta$ with $L_0R^2/\eta$ instead of $L_0R/\varepsilon$.\par
\end{center}

\endgroup
\end{document}